\documentclass{amsart}
\usepackage{amsthm,amssymb,amsmath,color}
\usepackage{pdfsync}

\newcommand{\levy}{L\'evy}

\newcommand{\clt}{central limit theorem}

\newcommand{\garch}{{\rm GARCH}$(1,1)$}

\newcommand{\ex}{{\rm e}\,}

\newcommand{\asy}{asymptotic}

\newcommand{\ts}{time series}
\newcommand{\tsa}{\ts\ analysis}

\definecolor{darkblue}{rgb}{.1, 0.1,.8}
\definecolor{darkgreen}{rgb}{0,0.8,0.2}
\definecolor{darkred}{rgb}{.8, .1,.1}

\newtheorem{lemma}{Lemma}[section]

\newtheorem{theorem}[lemma]{Theorem}

\newtheorem{proposition}[lemma]{Proposition}
\newtheorem{definition}[lemma]{Definition}
\newtheorem{corollary}[lemma]{Corollary}
\newtheorem{example}[lemma]{Example}
\newtheorem{exercise}[lemma]{Exercise}
\newtheorem{remark}[lemma]{Remark}
\newtheorem{fig}[lemma]{Figure}
\newtheorem{tab}[lemma]{Table}

\newcommand{\MC}{Markov chain}

\newcommand{\bfR}{{\bf R}}

\newcommand{\bth}{\begin{theorem}}
\newcommand{\ethe}{\end{theorem}}
\newcommand{\sv}{stochastic volatility}
\newcommand{\bre}{\begin{remark}\em }
\newcommand{\ere}{\end{remark}}

\newcommand{\ble}{\begin{lemma}}
\newcommand{\ele}{\end{lemma}}
\newcommand{\sre}{stochastic recurrence equation}
\newcommand{\pp}{point process}
\newcommand{\bde}{\begin{definition}}
\newcommand{\ede}{\end{definition}}
\newcommand{\bco}{\begin{corollary}}
\newcommand{\eco}{\end{corollary}}

\newcommand{\bpr}{\begin{proposition}}
\newcommand{\epr}{\end{proposition}}

\newcommand{\bexer}{\begin{exercise}}
\newcommand{\eexer}{\end{exercise}}

\newcommand{\bexam}{\begin{example}}
\newcommand{\eexam}{\end{example}}

\newcommand{\bfi}{\begin{fig}}
\newcommand{\efi}{\end{fig}}

\newcommand{\btab}{\begin{tab}}
\newcommand{\etab}{\end{tab}}

\newcommand{\lhs}{left-hand side}
\newcommand{\fidi}{finite-dimensional distribution}
\newcommand{\rv}{random variable}

\newcommand{\var}{{\rm var}}

\newcommand{\cov}{{\rm cov}}

\newcommand{\as}{{\rm a.s.}}

\newcommand{\rhs}{right-hand side}

\newcommand{\beao}{\begin{eqnarray*}}
\newcommand{\eeao}{\end{eqnarray*}\noindent}

\newcommand{\beam}{\begin{eqnarray}}
\newcommand{\eeam}{\end{eqnarray}\noindent}

\newcommand{\beqq}{\begin{equation}}
\newcommand{\eeqq}{\end{equation}\noindent}

\newcommand{\bce}{\begin{center}}
\newcommand{\ece}{\end{center}}

\newcommand{\barr}{\begin{array}}
\newcommand{\earr}{\end{array}}

\newcommand{\stp}{\stackrel{P}{\rightarrow}}

\newcommand{\stas}{\stackrel{\rm a.s.}{\rightarrow}}

\newcommand{\stv}{\stackrel{v}{\rightarrow}}

\newcommand{\vague}{\stackrel{\lower0.2ex\hbox{$\scriptscriptstyle
                    \it{v} $}}{\rightarrow}}
\newcommand{\weak}{\stackrel{\lower0.2ex\hbox{$\scriptscriptstyle
                    \it{w} $}}{\rightarrow}}
\newcommand{\what}{\stackrel{\lower0.2ex\hbox{$\scriptscriptstyle
                    \it{\hat{w}} $}}{\rightarrow}}

\newcommand{\bdis}{\begin{displaymath}}
\newcommand{\edis}{\end{displaymath}\noindent}

\newcommand{\N}{\mathbb{N}}
\newcommand{\R}{\mathbb{R}}

\newcommand{\nto}{n\to\infty}
\newcommand{\kto}{k\to\infty}
\newcommand{\xto}{x\to\infty}

\newcommand{\ov}{\overline}

\newcommand{\vep}{\varepsilon}

\newcommand{\la}{\lambda}

\newcommand{\regvary}{regularly varying}
\newcommand{\slvary}{slowly varying}
\newcommand{\regvar}{regular variation}

\newcommand{\bbr}{{\mathbb R}}

\newcommand{\bbz}{{\mathbb Z}}
\newcommand{\bbn}{{\mathbb N}}

\newcommand{\bbf}{{\mathcal F}}

\newcommand{\con}{convergence}

\newcommand{\st}{such that}

\newcommand{\fct}{function}

\newcommand{\ds}{distribution}

\newcommand{\rep}{representation}

\newcommand{\seq}{sequence}

\newcommand{\pro}{probabilit}

\newcommand{\ms}{measure}

\newcommand{\ld}{large deviation}

\newcommand{\bfX}{{\bf X}}
\newcommand{\bfM}{{\bf M}}

\newcommand{\bfy}{{\bf y}}

\def\1{\ensuremath{\mathrm{1}\hspace{-.35em} \mathrm{1}}} 

\def\E{{\mathbb E}}

\def\N{\mathbb{N}}
\def\P{{\mathbb{P}}}
\def\R{\mathbb{R}}

\renewcommand{\le}{\ensuremath{\leqslant}}
\renewcommand{\ge}{\ensuremath{\geqslant}}

\newcommand{\introo}[2]{{\left]{#1,\,#2\,}\right[\kern1pt}}

\newcommand{\intrfo}[2]{{\left[{#1,\,#2}\right[\kern1pt}}

\begin{document}
\today
\title[Precise large deviations for dependent regularly varying sequences ]
{Precise large deviations for dependent regularly varying sequences }
\thanks{Thomas Mikosch's research is partly supported
by the Danish Research Council (FNU) Grants 272-06-0442 and 09-072331. The
research of Thomas Mikosch and Olivier Wintenberger is partly supported by a
Danish-French Scientific Collaboration Grant of the French Embassy in
Denmark. Both authors would like to thank their home institutions for
hospitality when visiting each other.}
\author[T. Mikosch]{Thomas Mikosch}
\author[O. Wintenberger]{Olivier Wintenberger}
\address{Thomas Mikosch, University of Copenhagen, Department of Mathematics,
Universitetsparken 5,
DK-2100 Copenhagen\\ Denmark} \email{mikosch@math.ku.dk}

\address{Olivier Wintenberger,
            Universit\'e de Paris-Dauphine,
            Centre De Recherche en Math\'ematiques de la D\'ecision
            UMR CNRS 7534,
            Place du Mar\'echal De Lattre De Tassigny,
            75775 Paris Cedex 16, France}
            \email{owintenb@ceremade.dauphine.fr}

\maketitle
\begin{abstract}
We study a precise \ld\ principle for a stationary \regvary\ \seq\
of \rv s. This principle extends the classical results of
A.V. Nagaev \cite{nagaev:1969} and S.V. Nagaev \cite{nagaev:1979}
for iid \regvary\ \seq s. The proof uses an idea of
Jakubowski \cite{jakubowski:1993,jakubowski:1997} in the context of
\clt s with infinite variance stable limits.
We illustrate the principle for \sv\ models, functions of a \MC\ 
satisfying a polynomial drift condition and solutions of linear
and non-linear \sre s.
\end{abstract}



\section{Introduction}
The aim of this paper is to study {\em precise \ld\ \pro ies}
for \seq s of dependent and heavy-tailed \rv s. To make the notion of
heavy tails  precise, we assume that  the stationary \seq\ $(X_t)$
has \regvary\ \fidi s in the sense defined in Section~\ref{subsec:rv}.
A particular con\seq\ is that the \ds\ of a generic variable $X$
of this
\seq\ has \regvary\ tails. This means that there exist $\alpha>0$,
$p,q\ge 0$ with $p+q=1$ and a \slvary\ \fct\ $L$ \st
\beam\label{eq:tailbalance}
\dfrac{\P(X>x)}{\P(|X|>x)}\sim p\,\dfrac{L(x)}{x^\alpha}\quad
\mbox{and}
\quad \dfrac{\P(X\le -x)}{\P(|X|>x)}\sim
q\,\dfrac{L(x)}{x^\alpha}\,,\quad \xto .
\eeam
The latter condition is often referred to as a {\em tail balance
  condition.}
\par
In the case of an iid \seq\ satisfying \eqref{eq:tailbalance}
one can derive precise \asy\ bounds for the tails of the
random walk $(S_n)$ with step \seq\  $(X_t)$ given by
\beao
S_0=0\quad\mbox{and} \quad S_n=X_1+\cdots+X_n\,,\quad n\ge 1\,.
\eeao
We recall a classical result which can be found in the papers of
A.V. and S.V. Nagaev \cite{nagaev:1969,nagaev:1979} and
Cline and Hsing \cite{cline:hsing:1998}.
\bth\label{thm:nagaev}
Assume that $(X_i)$ is an iid \seq\ with a \regvary\ \ds \
in the sense of \eqref{eq:tailbalance}.
Then the following relations
hold for $\alpha>1$ and suitable \seq s $b_n\uparrow\infty$:
\beam\label{eq:lda}
\lim_{\nto}
\sup_{x\ge b_n}
\left|\dfrac{\P (S_n- \E S_n>x)}{n\,\P(|X|>x)}-
p\right|=0
\eeam
and
\beam\label{eq:ldb}
\lim_{\nto}
\sup_{x\ge b_n}\left|\dfrac{\P(S_n- \E S_n\le -x)}{n\,\P(|X|>x)}-
q\right|=0\,.
\eeam
If $\alpha>2$ one can choose $b_n=\sqrt{a n \log n}$, where
$a>\alpha-2$, and for $\alpha\in (1,2]$, $b_n=n^{\delta+1/\alpha}$
for any $\delta>0$. For $\alpha\le 1$, \eqref{eq:lda} and
\eqref{eq:ldb} remain valid with $\E S_n$ replaced by $0$
and one can choose $b_n=n^{\delta+1/\alpha}$ for any $\delta>0$.
\ethe
We call results of the type \eqref{eq:lda} and
\eqref{eq:ldb} a {\em precise \ld\  principle} in contrast to the
majority of results in \ld\ theory where the logarithmic \pro ies
$n^{-1}\log \P(n^{-1}(Y_n- \E Y_n) \in A)$ are studied for sets $A$ bounded
away from zero and
suitable \seq s
$(Y_n)$ of \rv s (not necessarily constituting a random walk)
or even random elements taking values in some
abstract spaces; see e.g. the monograph by Dembo and Zeitouni \cite{dembo:zeitouni:2010}.
As a matter of fact, precise \ld\ principles can be derived for
iid heavy-tailed \seq s more general than \regvary\ ones, e.g. for the
general class of random walks $(S_n)$ with subexponential steps; see
e.g. Cline and Hsing \cite{cline:hsing:1998},
Denisov et al. \cite{denisov:dieker:shneer:2008}, Mogulskii
\cite{mogulski:2009} and the references cited therein. We also mention that Theorem~\ref{thm:nagaev}
can be extended to iid \regvary\ random vectors (see Section~\ref{subsec:rv}
for a definition) and an analog of Donsker's theorem for large
deviations in Skorokhod space can be
proved as well; see Hult et al. \cite{hult:lindskog:mikosch:samorodnitsky:2005}.
\par
Theorem~\ref{thm:nagaev} serves as a benchmark result for the
purposes of this paper.
In this paper we extend
Theorem~\ref{thm:nagaev} to suitable \regvary\ stationary \seq s
$(X_t)$.
Various examples of precise \ld\ principles
have been derived in the literature. Under rather general dependence
conditions on the \regvary\ \seq\ $(X_t)$ with index $\alpha<2$, Davis and Hsing
\cite{davis:hsing:1995} and Jakubowski \cite{jakubowski:1993,jakubowski:1997}
proved the existence of some \seq s $(b_n)$ \st\ $b_n^{-1}S_n\stp 0$
and
\beam\label{eq:kkk}
\lim_{\nto}\dfrac{\P(S_n >b_n)}{n\,\P(|X|>b_n)}\,.
\eeam
They could in general not specify the order of magnitude of the
\seq s $(b_n)$. The method of
proof for these results
could not be extended to the case $\alpha\ge 2$. Moreover, work of
Lesigne and Voln\'y \cite{lesigne:volny:2001} indicates that results
of the type of Theorem~\ref{thm:nagaev} may fail for certain stationary ergodic
martingale difference \seq s. To be more precise, they proved that
$\limsup_{\nto}\P(S_n >n)/[n\,\P(|X|>n)]=\infty$ is possible for
such \seq s. Gantert \cite{gantert:2000} proved \ld\ results of
logarithmic type for stationary ergodic \seq s $(X_t)$ satisfying a geometric
$\beta$-mixing condition. The latter condition ensures that the tail \asy s do
not differ from the iid case.
\par
An analog of Theorem~\ref{thm:nagaev} for linear processes
$X_t=\sum_{j=0}^\infty \psi_j Z_{t-j}$, $t\in\bbz$,
under suitable assumptions on
the \seq\ of  real numbers $(\psi_j)$ (ensuring the existence of the
infinite series) and assuming \regvar\ of the iid
innovations $(Z_t)$
was proved in Mikosch and Samorodnitsky
\cite{mikosch:samorodnitsky:2000}. The limiting constants $p$ and $q$ in \eqref{eq:lda} and
\eqref{eq:ldb}, respectively,  had to be replaced by quantities depending on $p,q$ and
the \seq\ $(\psi_j)$. The region $(b_n,\infty)$, where the
\ld\ principle holds, remains the same as for an iid
\regvary\ \seq .
\par
Similar results were obtained in Konstantinides and
Mikosch \cite{konstantinides:mikosch:2005} for solutions to the \sre\
$X_t=A_tX_{t-1}+B_t$, $t\in\bbz$,  with  iid $((A_t,B_t))_{t\in\bbz}$
with a generic element $(A,B)$, $A,B\ge 0$ a.s.,
$B$ \regvary\ with index $\alpha>0$ and $EA^\alpha<1$. They showed that the limits \eqref{eq:kkk} exist
and are positive  for \seq s $(b_n)$ comparable to those in Theorem~\ref{thm:nagaev}; uniform results like in \eqref{eq:lda} and
\eqref{eq:ldb} were not achieved. For the same type of \sre\ with $B$
not necessarily positive, Buraczewski et
al. \cite{buraczewski:damek:mikosch:zienkiewicz:2011} proved precise \ld\
principles. The main difference to
\cite{konstantinides:mikosch:2005} is the assumption that $(X_t)$ is
\regvary\ with some positive index $\alpha$ while
$(A_t,B_t)$ has moments of order $\alpha+\delta$ for some positive
$\delta$. In this case, the celebrated paper of Kesten
\cite{kesten:1973}, under appropriate conditions on the \ds\ of
$(A,B)$, yields  that $(X_t)$ is indeed \regvary\ with index $\alpha$; see also
Goldie \cite{goldie:1991}. It is shown in
\cite{buraczewski:damek:mikosch:zienkiewicz:2011}
that the relation
\beao
\limsup_{\nto}\sup_{x\ge b_n}\dfrac{\P(S_n >x)}{n\,\P(|X|>x)}<\infty
\eeao
holds for suitable \seq s $b_n\to\infty$ \st\ $b_n^{-1}S_n\stp
0$. Again, the \seq s $(b_n)$ are close to those in Theorem~\ref{thm:nagaev}.
However, uniform relations of type \eqref{eq:lda} and \eqref{eq:ldb} are not
true in the unbounded regions $(b_n,\infty)$ but in bounded regions
$(b_n,c_n)$ \st\ $b_n\to\infty$ and $c_n=\ex^{s_n}$ for $s_n\to\infty$
and $s_n=o(n)$.
\par
In this paper, we will approach the problem of precise \ld s
from a more general point of view. A key idea for this approach
can be found in the papers of Jakubowski
\cite{jakubowski:1993,jakubowski:1997}, where this idea was used to
prove central limit theory with infinite variance stable limits for the
partial sums $(S_n)$ of a general stationary \seq\ with \regvary\
marginals; see also the recent paper Bartkiewicz et
al. \cite{bartkiewicz:jakubowski:mikosch:wintenberger:2011}, where the
same idea was exploited.
The following inequality is crucial for proving the results of
this paper: for every $k\ge 2$, some constant $b_+$,
\beam\label{eq:gg}
\lefteqn{\Big|\dfrac{\P(S_n>x)}{n\,\P(|X|>x)}-b_+\Big|}\nonumber\\
&\le &\Big|
\dfrac{\P(S_n> x)-n\,(\P(S_{k+1}> x)-\P(S_k> x))}{n\,\P(|X|> x)}\Big|
+\Big|\dfrac{\P(S_{k+1}> x)-\P(S_k> x)}{\P(|X|>
  x)}-b_+\Big|\,.
\eeam
Regular variation of $(X_t)$ ensures that the second
quantity in \eqref{eq:gg} is negligible, by first letting $\xto$ and
then $\kto$. The first expression in \eqref{eq:gg} provides a link
between the \asy s of the tail  $\P(S_n>x)$ for increasing values of
$n$, $x\ge b_n$ and
the \regvary\ tails $\P(S_k>x)$ and $\P(S_{k+1}>x)$ for every fixed
$k$. Thus the tail \asy s of $\P(S_n>x)$ are derived from the known tail
\asy s for finite sums, again by first letting $\nto$ and then $\kto$.
\par
This paper is organized as follows. In Section~\ref{subsec:prel}
we introduce some of the basic conditions and notions needed
throughout the paper.
These include \regvar\ of a stationary \seq\ and
an anti-clustering condition. In Section~\ref{sec:main} we prove the
main result of this paper: Theorem~\ref{th:main} provides a general precise
\ld\ principle for \regvary\ stationary \seq s.
 Under \regvar\ and anti-clustering
conditions we will show precise \ld\ principles of the
following type:
\beam\label{eq:1}
\lim_{n\to \infty}\sup_{x\in \Lambda_n}\Big|\frac{\P(S_n> x)}{n\,\P(|X|> x)}-b_+\Big|=0\,,
\eeam
for some non-negative constant $b_+$ and a \seq\ of sets  $\Lambda_n
\subset (0,\infty)$ \st\ $b_n=\inf \Lambda_n\to\infty$.  In Section~\ref{sec:exam}
we will apply the \ld\ principle \eqref{eq:1} to a variety of
important \regvary\ \ts\ models, including the \sv\ model, solutions
to \sre s and functions of Markov chains. These are examples of rather different 
dependence structures, showing that the \ld\ principle does not depend on a 
particular mixing condition or on the Markov property. 
\par
However, we give 
special emphasis to functions of a Markov chain satisfying a polynomial drift condition. 
Theorems \ref{thm:mc} and \ref{thm:int}
are our main results for \MC s.  Theorem \ref{thm:mc} is obtained 
by a direct application of Theorem~\ref{th:main}, exploiting a sophisticated
exponential bound for partial sums of \MC s due to Bertail and Cl\'emencon
\cite{bertail:clemencon:2009}. Theorem \ref{thm:mc} implies 
Theorem \ref{thm:int}. It 
yields an intuitive interpretation of relation
\eqref{eq:1} in terms of the regeneration property of $(X_t)_{t=1,\ldots,n}$. 
Given an atom $A$ of the underlying chain, one can split the chain into
a random number $N_A(n)$ of iid random cycles. Denoting the block sum of the 
$X_t$'s over the $i$th cycle by $S_{A,i}$, it will be shown that the iid 
$S_{A,i}$'s inherit \regvar\ from $X$, and then we can apply the classical
result of Theorem~\ref{thm:nagaev} to 
$\P\Big(\sum_{i=1}^{N_A(n)-1} S_{A,i}>x\Big)$. If $b_+>0$ the tails $\P_A(S_{A,1}>x)$ and $\P(|X|> x)$ are equivalent. There is a major difference between 
an iid \seq\ and the dependent \seq\ $(X_t)$: if the first generation 
time $\tau_A$ is larger than $n$, it has significant influence on the region 
$\Lambda_n$, where \eqref{eq:1} holds. It turns out that one has 
for any $x\ge b_n$,
\beao
\dfrac{\P(S_n>x)}{n\P(|X|>x)}\sim b_+ +\dfrac{\P(S_n>x,\tau_A>n)}{n\P(|X|>x)}\,,
\eeao
and the second term is in general not negligible, leading to the fact that 
\eqref{eq:1} may only be valid in a bounded  region $(b_n,c_n)$. Thus we found
an explanation for the same observation we experienced in the case 
of a \MC\
given by a \sre ; see the discussion above.

\section{Preliminaries}\label{subsec:prel}
\setcounter{equation}{0}
\subsection{Regular variation}\label{subsec:rv}
Throughout this paper we assume that $(X_t)$ is stationary. Such a
\seq\ is
{\em \regvary\ with index $\alpha>0$} if
the \fidi s
of $(X_t)$ have a jointly \regvary\ \ds\ in the following sense: for
every $d\ge 1$, there exists a non-null Radon \ms\ $\mu_d$ on the
Borel $\sigma$-field of $\ov \bbr^d\backslash\{\bf0\}$, where
$\ov \bbr=\bbr\cup\{\pm \infty\}$,
(this means
that $\mu_d$ is finite on sets bounded away from zero) \st\
\beao
n\,\P(a_n^{-1}(X_1,\ldots,X_d)\in \cdot )\stv \mu_d(\cdot)\,, \eeao
where $\stv$ denotes vague \con\ (see  e.g.
\cite{kallenberg:1983,resnick:1987}) and $(a_n)$
satisfies
$n\,\P(|X|>a_n)\sim 1$.
The limiting \ms s have the property $\mu_d(tA)=t^{-\alpha}\mu_d(A)$,
$t>0$, for any Borel set $A$. We refer to $\alpha$ as the {\em index of
\regvar } of $(X_t)$ and its \fidi s. We refer to Basrak and Segers
\cite{basrak:segers:2009} for an insightful description of \regvar\ for stationary processes.
\par
In what follows, we refer to condition ${\bf RV}_\alpha$ if $(X_t)$
satisfies the conditions above for some $\alpha>0$ and a \seq\ of
limiting \ms s $(\mu_d)$.
\par
In Section~\ref{sec:exam} we will consider some prominent examples of
\regvary\ \ts .
\par
The \regvar\ property of $(X_t)$ implies that
the limits
\beam\label{eq:b+}
b_+(k)=\lim_{\xto}\dfrac{\P(S_k> x)}{\P(|X|> x)}
=\lim_{\nto}n\,\P(S_k> a_n)
,\qquad k\ge 1,
\eeam
exist. These quantities play a crucial role in our investigations on
\ld s;  see for example Theorem~\ref{th:main}. The
limiting constants
\beam\label{eq:b-}
b_-(k)=\lim_{\xto}\dfrac{\P(S_k\le -x)}{\P(|X|> x)}
=\lim_{\nto}n\,\P(S_k\le -a_n)
,\qquad k\ge 1,
\eeam
also exist by virtue of \regvar\ of $(X_t)$.
\par
In our main result Theorem~\ref{th:main} we require that the
limit
\beao
b_+=\lim_{\kto} (b_+(k+1)-b_+(k))
\eeao
exists; the existence of $b_+$ does not directly follow from \regvar\ of
$(X_t)$. In the examples of Section~\ref{sec:exam} we show that $b_+$
is easily calculated for some major \ts\ models.
If $b_+$ exists it is non-negative since it is the limit of a
C\`esaro mean: $b_+=\lim_{\kto} k^{-1} b_+(k)$.
\par
The constants $b_+$ and $b_-$ (the latter constant is defined in the straightforward
way) figure prominently in \asy\ results
for the partial sums $(S_n)$ with infinite variance stable
limits. Indeed, the \levy\ \ms\ $\nu$ of the stable limit has \rep\
$\nu(x,\infty)=b_+x^{-\alpha}$ and $\nu(-\infty,-x)=b_-x^{-\alpha}$,
$x>0$; see Bartkiewicz et al.~\cite{bartkiewicz:jakubowski:mikosch:wintenberger:2011}.

\subsection{Anti-clustering condition}\label{subsec:ac}
Assume that $(X_t)$ satisfies the \regvar\ condition ${\bf
  RV}_\alpha$. For studying the limit theory for the extremes of
dependent \seq s it is common to assume {\em anti-clustering conditions}; see
e.g. Leadbetter et al.
\cite{leadbetter:lindgren:rootzen:1983}, Leadbetter and Rootz\' en
\cite{leadbetter:rootzen:1988} and Embrechts et al.
\cite{embrechts:kluppelberg:mikosch:1997}, Chapter 5.
These conditions ensure that possible clusters of exceedances of high thresholds by
the \seq\ $(X_t)$ cannot be too large. In other words, ``long-range
dependencies of extremes'' are avoided. Anti-clustering conditions are
also needed for proving \asy\ theory for partial sums with infinite variance stable
limits; see Davis and Hsing
\cite{davis:hsing:1995}, Jakubowski
\cite{jakubowski:1993,jakubowski:1997}, Basrak
and Segers \cite{basrak:segers:2011}, and
Bartkiewicz et
al.~\cite{bartkiewicz:jakubowski:mikosch:wintenberger:2011}. In the
latter reference, the different conditions are discussed and compared.
Davis and Hsing \cite{davis:hsing:1995} and Jakubowski
\cite{jakubowski:1993,jakubowski:1997} also proved \ld\ results in the
case $\alpha<2$ under anti-clustering conditions.
\par
We introduce the
following anti-clustering condition which is close to those
in the literature mentioned above.\\[1mm]
{\em Condition} ${\bf AC}_\alpha$: There exist
$\delta_k\downarrow 0$ as $k\to\infty$ and a \seq\ of sets
$\Lambda_n\subset (0,\infty)$, $n=1,2,\ldots$, with $b_n=\inf \Lambda_n$
\st\ $n\,\P(|X|>b_n)\to 0$  as $\nto$ and
\beao
\lim_{k\to\infty}\limsup_{\nto} \sup_{x\in \Lambda_n}
\delta_k^{-\alpha}\sum_{j= k}^n\P(|X_j|> x\delta_k\mid |X_0|> x\delta_k)=0\,.
\eeao
This condition is tailored for the purposes of our paper: the sets
$(\Lambda_n)$ with $\lim_{\nto}b_n=\infty$ are those
which appear in the precise \ld\ results \eqref{eq:1}.
\par
Condition  ${\bf AC}_\alpha$  is easily verified for the examples of \ts\
models in Section~\ref{sec:exam}.

\section{Main result}\label{sec:main}
\setcounter{equation}{0}
In this section we formulate and prove the main result on precise \ld\
principles for \regvary\ stationary \seq s.
\bth\label{th:main} Assume that the stationary \seq\ $(X_t)$ of
real-valued \rv s satisfies the following  conditions.
\begin{enumerate}
\item
The \regvar\ condition ${\bf RV}_\alpha$  for some $\alpha>0$.
\item
The anti-clustering condition ${\bf AC}_\alpha$ for a \seq\
$\delta_k=o(k^{-2})$, $\kto$, and sets $(\Lambda_n)$ \st\ $b_n=\inf
\Lambda_n\to \infty$ as $\nto$.
\item
The limit
$
b_+=\lim_{k\to\infty}(b_+(k+1)-b_+(k))
$ exists, where the constants $(b_+(k))$ are defined in \eqref{eq:b+}.
\item
For the \seq s  $(\Lambda_n)$, $(\delta_k)$ from ${\bf AC}_\alpha$ and a
\seq\ $(\vep_k)$ satisfying
$\varepsilon_k=o(k^{-1})$ and $(k+1)\delta_k\le \vep_k$,
\beam\label{eq:main}
\lim_{\kto}\limsup_{\nto}\sup_{x\in \Lambda_n}
\frac{\P\Big(\sum_{i=1}^n X_i\1_{\{|X_i|\le \delta_k \,x\}}>\varepsilon_k x\Big)}{n\,\P(|X|>x)}=0.
\eeam
\end{enumerate}
Then the large deviation principle \eqref{eq:1} holds.
\ethe
The corresponding result for the left tails $\P(S_n\le -x)$, $x>0$, is obtained
by replacing the variables $X_t$ by $-X_t$, $t\in\bbz$. Then one also needs to
assume that the limit $b_-$ exists which is defined correspondingly.
\bre\label{rem:less1}
In the case $\alpha<1$, \eqref{eq:main} is satisfied for suitable choices of
$(\delta_k)$ and $(\vep_k)$.
Indeed, an application of Markov's inequality and
Karamata's theorem
(see Bingham at al. \cite{bingham:goldie:teugels:1987})
yields uniformly for $x>b_n$,
\beao
\P\Big(\sum_{i=1}^n X_i\1_{\{|X_i|\le \delta_k \,x\}}>\varepsilon_k
x\Big)
&\le & (x\vep_k)^{-1}\,n\,\E |X|\1_{\{|X|\le \delta_k \,x\}}\\
&\sim &\delta_k^{1-\alpha}\vep_k^{-1}\,n\,\P(|X|>x)\,.
\eeao
Thus \eqref{eq:main} is satisfied for $\Lambda_n=(b_n,\infty)$ if we
choose e.g. $\delta_k=\ex^{-k}$ and $\vep_k=k^{-2}$.
\ere
\bre\label{rem:garch} Assume $\alpha\in (0,2)$ and
$(X_t)$ conditionally independent and symmetric
given some $\sigma$-field $\bbf$. This condition is often
satisfied in models of the type
$X_t=\sigma_t\,Z_t$ with iid symmetric $(Z_t)$, for example
if $(Z_t)$ and $(\sigma_t)$ are
independent; see the \sv\ model of
Section~\ref{subsec:sv}.
Alternatively, if $(\sigma_t)$ is
predictable
with respect to the filtration generated by the \seq\ $(Z_t)$ then
$(X_t)$ is conditionally independent and symmetric. Prominent examples
of this type are GARCH-type models, where
$(Z_t)$ is often assumed iid standard normal or student distributed.
Indeed, first applying the Chebyshev inequality conditional on
$\bbf$ and then taking expectations, we obtain by Karamata's theorem
(see Bingham at al. \cite{bingham:goldie:teugels:1987})
uniformly for $x\in \Lambda_n=(b_n,\infty)$,
\beao
\P\Big(\sum_{i=1}^n X_i\1_{\{|X_i|\le \delta_k \,x\}}>\varepsilon_k
x\Big)&\le& (\vep_k x) ^{-2 }n\, \E X^2\1_{\{|X|\le \delta_k \,x\}}\\
&\sim & \delta_k^{2-\alpha}\vep_k^{-2 }n\, \P (|X|> x)\,.
\eeao
Thus \eqref{eq:main} holds e.g. for $\delta_k=\ex^{-k}$ and $\vep_k=k^{-2}$.
\ere
\bre\label{rem:0} Recall that $(b_n)$ is chosen \st\
$n\,\P(|X|>b_n)\to 0$. For an  iid $(X_t)$, this condition is
necessary for the weak law of large numbers $b_n^{-1}S_n\stp 0$.
Under this and some other mild conditions,
we may assume without loss of generality
that the  \rv s  $(X_i\1_{\{|X_i|\le \delta_k \,x\}})$ in
\eqref{eq:main} are mean corrected.
Indeed, we will prove that
\beam\label{eq:mild}
n\,\sup_{x\in\Lambda_n} x^{-1}|\E X \1_{\{|X|\le x\}}|=o(1)\,,\quad\nto\,.
\eeam
This condition is trivial if $X$ is symmetric.\\
{\em The case $\alpha<1$.}
By Karamata's theorem
and the choice of $(b_n)$,
\beao
n\,|\E X\1_{\{|X|\le x\}}|\le n\,\E |X|\1_{\{|X|\le x\}}  \sim  c\,n\,x\, \P(|X|>x)\le c\,x \,n\P(|X|>b_n)=o(x)\,.
\eeao
{\em Here and in what follows, we write $c$ for any positive constants
 whose value is not of interest, for example, the same $c$ may denote different
 constants in the same formula.}\\
{\em The case $\alpha=1$.} If $\E X=0$ and $n=O(b_n)$ then
$n\,|\E X\1_{\{|X|\le x\}}|=o(n)=o(x)$. If $\E |X|=\infty$,
$E|X|\1_{\{|X|\le x\}}$
is a \slvary\ \fct , and therefore for large $n$ and any small $\epsilon>0$,
$n\,|EX\1_{\{|X|\le x\}}|\le n\,x^\epsilon$. If $b_n=n^{1+\delta}$ for
some $\delta>0$, choosing $\epsilon$ sufficiently small, we obtain
$n\,|EX\1_{\{|X|\le x\}}|=o(x)$.\\
{\em The case $\alpha>1$.}
By Karamata's theorem, since $EX=0$
and by the
choice of $(b_n)$, as $\nto$,
\beao
n\,|\E X\1_{\{|X|\le x\}}|&=&
n\,|\E X\1_{\{|X|> x\}}|\le n\,\E|X|\1_{\{|X|> x\}}\\
&\sim& c\,
n\,x\,\P(|X|>x)
\le c\,x\,[n\P(|X|>b_n)]=o(x)\,.
\eeao
\ere

\begin{proof}
We have for fixed $k\ge 2$,
\beao\lefteqn{
\sup_{x\in \Lambda_n}
\Big|\dfrac{\P(S_n>x)}{n\,P(|X|>x)}-b_+\Big|}\\
&\le &\sup_{x\in\Lambda_n}\Big|
\dfrac{\P(S_n> x)-n\,(\P(S_{k+1}> x)-\P(S_k> x))}{n\,\P(|X|> x)}\Big|
+\sup_{x\in\Lambda_n}\Big|\dfrac{\P(S_{k+1}> x)-\P(S_k> x)}{\P(|X|>
  x)}-b_+\Big|\\
&=& I_{1,k}+I_{2,k}\,.
\eeao
By \regvar\ of $(X_t)$, the limit
\beao
\lim_{n\to\infty} I_{2,k}= |(b_{+}(k+1)-b_{+}(k))-b_+|
\eeao
exists for every $k\ge 2$ and any \seq\ $(\Lambda_n)$ \st\
$\inf\Lambda_n\to\infty$. By assumption,
$\lim_{\kto} |(b_{+}(k+1)-b_{+}(k)-b_+|=0$. Therefore it suffices to
study the \asy\ behavior of $I_{1,k}$.
\par
For any $\delta>0$ and $x>0$, consider
\beao
\overline X_i=X_i\1_{\{|X_i|\le  x\delta\}}
\quad\mbox{and}\quad
\underline X_i=X_i\1_{\{|X_i|>  x\delta\}}\,,\quad i=1,2,\ldots.
\eeao
and for $n\ge 1$,
\beao
\ov S_n=\sum_{i=1}^n \overline X_i\quad\mbox{and}\quad
\underline S_n=\sum_{i=1}^n \underline X_i\,.
\eeao
Then, for any
$\varepsilon\in (0,1)$ and $j\ge 1$,
\beao
\P(\underline S_j> (1+\varepsilon)\,x)-
\P(-\overline S_j> \varepsilon \,x)\le
\P(S_j> x)\le \P(\underline S_j> (1-\varepsilon)\,x)+
\P(\overline S_j> \varepsilon\, x).
\eeao
Multiple application of these inequalities yields
\beao
A_1+A_2+A_3\le \dfrac{\P(S_n> x)-n\,(\P(S_{k+1}> x)-\P(S_k> x))}{n\,\P(|X|> x)}
\le B_1+B_2+B_3\,,
\eeao
where
\begin{eqnarray*}
A_1&=&\frac{\P(\underline S_n> (1+\varepsilon)x)-n\,( \P(\underline
  S_{k+1}>  (1+\varepsilon)x)-\P(\underline S_k>
  (1+\varepsilon)x))}{n\,\P(|X|> x)}\,,\\
A_2&=&\frac{-\P(-\overline S_n> \varepsilon x)-n\,( \P(\overline
  S_{k+1}>  \varepsilon x)-\P(-\overline S_k>   \varepsilon
  x))}{n\,\P(|X|> x)}\,,\\
A_3&=&\frac{\P(\underline S_{k+1}>  (1+\varepsilon) x)-
\P(\underline S_{k+1}>  (1-\varepsilon) x)}{\P(|X|> x)}\,,\\
B_1&=&\frac{\P(\underline S_n> (1-\varepsilon)x)-n\,( \P(\underline S_{k+1}>  (1-\varepsilon)x)-\P(\underline S_k>  (1-\varepsilon)x))}{n\,\P(|X|> x)}\,,\\
B_2&=&\frac{\P(\overline S_n> \varepsilon x)+n\,( \P(-\overline S_{k+1}>  \varepsilon x)+\P(\overline S_k>   \varepsilon x))}{n\,\P(|X|> x)}\,,\\
B_3&=&\frac{\P(\underline S_{k+1}>  (1-\varepsilon) x)-\P(\underline S_{k+1}>  (1+\varepsilon) x)}{\P(|X|> x)}.
\end{eqnarray*}
We will derive upper bounds for the $B_i$'s. Lower bounds for the
$A_i$'s can be derived in the same way and are therefore omitted.
\par
An application of Jakubowski
\cite{jakubowski:1997}, Lemma~3.2, to the stationary \seq\ $(\underline X_t)$
yields
for fixed $k\ge 2$, $x,\delta,\vep>0$,
\beao
|B_1|&\le& 3\,\frac{k\,\P(|X|>\delta \,x )}{n\,\P(|X|> x )}+2\,\sum_{j=
  k}^n\frac{\P(|X_j|>\delta\,x ,|X_0|>\delta\,x)}{\P(|X|> x)}\\
&=& B_{11}+B_{12}\,.
\eeao
In view of \regvar\ of $X$, $\P(|X|>\delta \,x )/\P(|X|>x)\to\delta^{-\alpha}$.
Hence
\beao
\limsup_{\nto }\sup_{x\in\Lambda_n}B_{11}=0\,,\quad k\ge 2\,,
\eeao
An application of ${\bf AC}_\alpha$ with $\delta=\delta_k$ yields that
$\lim_{\kto}\limsup_{\nto} \sup_{x\in\Lambda_n} B_{12}=0$.
Hence
\beao
\lim_{\kto}\limsup_{\nto} \sup_{x\in\Lambda_n} B_{1}=0\,.
\eeao
\par
Next consider $B_2$. In addition to the condition
$\vep=\vep_k=o(k^{-1})$ assume that
$(k+1)\delta_k\le \vep_k$. This choice is always possible since
we also assume $\delta=\delta_k=o(k^{-2})$.
Then
$|\overline S_{k+1}|\le \vep\,x$,
$\P(-\overline S_{k+1}>\varepsilon x)=\P(\overline S_{k}>\varepsilon
x)=0$ and $B_2$ degenerates to the expression
$
\P(\overline S_n>\varepsilon x)/(n\,\P(|X|>x))$. By assumption
\eqref{eq:main},
this condition is \asy ally negligible.
\par
Finally, consider $B_3$. Fix $k\ge 2$. In what follows, the constants
$\vep,\delta\in (0,1)$ will also depend on $k$.
Consider the sets
$$
A_{\gamma,\delta}(k)=\Big\{ \bfy\in \R^{k}:
\sum_{i=1}^k y_i\1_{\{|y_i| >\delta\}}>\gamma\Big\},\quad\gamma\,,\delta>0\,.
$$
Observe that
\beao
\{\underline S_{k+1}>  \gamma x\} =\{x^{-1}(X_1,\ldots,X_{k+1})\in
A_{\gamma,\delta}(k+1)\}\,,
\eeao
the sets $A_{\gamma,\delta}(k)$
are bounded away from $0$ and
$A_{\gamma,\delta}(k)=\gamma A_{1,\delta/\gamma}(k)$.
Condition ${\bf RV}_\alpha$ ensures the existence of the limit
\beao
&&\lim_{\xto}
B_3\\&=&\mu_{k+1}(A_{1-\varepsilon,\delta})-\mu_{k+1}(A_{1+\varepsilon,\delta})\\
&=&(1-\varepsilon)^{-\alpha}\mu_{k+1}(A_{1,\delta/(1-\varepsilon)})
-(1+\varepsilon)^{-\alpha}\mu_{k+1}(A_{1,\delta/(1+\varepsilon)})\\
&=&((1-\varepsilon)^{-\alpha} -(1+\varepsilon)^{-\alpha})
\mu_{k+1}(A_{1,\delta/(1-\varepsilon)})
-(1+\varepsilon)^{-\alpha}(\mu_{k+1}(A_{1,\delta/(1+\varepsilon)})-\mu_{k+1}(A_{1,\delta/(1-\varepsilon)}))\\
&=& B_{31}+B_{32}\,.
\eeao
By a Taylor expansion,
$B_{31}\le c\,\vep \,
\mu_{k+1}(A_{1,\delta/(1+\varepsilon)})$.
We observe that
\beam\label{eq:3}\lefteqn{
\Big\{\bfy\in \bbr^{k+1}:
\sum_{i=1}^ky_i>1+k\delta/ (1+\vep)\Big\}}\nonumber\\
&\subset&
A_{1,\delta/(1+\varepsilon)}=
\Big\{\bfy\in \bbr^{k+1}:
\sum_{i=1}^ky_i>1+\sum_{i=1}^ky_i\1_{\{|y_i|\le \delta/(1+\vep)\}}\Big\}\\
&\subset &\Big\{\bfy\in \bbr^{k+1}:
\sum_{i=1}^ky_i>1-k\delta/ (1+\vep)\Big\}\,.\nonumber
\eeam
Assume that $\delta=\delta_k=o(k^{-1})$ as $k\to\infty$. Then  for $k$
sufficiently large,
\beao
B_{31}\le c\,\vep \,
(1-k \delta/(1+\vep))^{-\alpha}\,\mu_{k+1}\Big(\Big\{\bfy\in
\bbr^{k+1}:\sum_{i=1}^{k+1}
y_i>1\Big\}\Big)\le  c\,\vep \, b_+(k+1)\,.
\eeao
Since we assume that $b_+$ exists a C\`esaro limit argument yields that
$\lim_{k\to\infty} k^{-1} b_+(k+1)=b_+$. Now choose
$\vep=\vep_k=o(k^{-1})$.
Then
$\lim_{k\to\infty}B_{31}=0.$
Similar arguments, using \eqref{eq:3}, yield
\beao
B_{32}&\le& c
\,b_+(k+1)\big((1-k\delta/(1-\vep))^{-\alpha}-
(1+k\delta/(1+\vep))^{-\alpha}\big)\\&\le&
c k\,\delta\, b_+(k+1) =o(1)\,,\quad\kto\,,
\eeao
provided $\delta=\delta_k=o(k^{-2})$. Thus we proved that
\beao
\lim_{\kto}\limsup_{\nto} \sup_{x\in\Lambda_n}B_3=0\,.
\eeao
\par

This concludes the proof.
\end{proof}

\section{Examples}\label{sec:exam}
In this section we want to apply Theorem~\ref{th:main} to a variety of
\ts\ models. Since there exists a calculus for multivariate
\regvar\ (e.g. Resnick  \cite{resnick:1987,resnick:2007}, Hult and
Lindskog \cite{hult:lindskog:2005,hult:lindskog:2006}, Basrak and
Segers \cite{basrak:segers:2009}) it is not difficult to
show the \regvar\ condition ${\bf RV}_\alpha$, the anti-clustering condition
${\bf AC}_\alpha$ and the existence of the limit $b_+=\lim_{\kto}(b_+(k+1)-b_+(k))$
for the examples below. However, it can take some efforts to prove
condition \eqref{eq:main}. In the iid case, one would use exponential
inequalities of Nagaev-Fuk or Prokhorov type; see e.g. the monograph
Petrov \cite{petrov:1995} for an overview of such inequalities. In the
case of dependent \seq s $(X_t)$ analogs of these inequalities exist,
but their application is not always straightforward; see e.g. the case
of \MC s in Section~\ref{subsec:mc} below.

\setcounter{equation}{0}
\subsection{$m_0$-dependent sequences}
In this section we consider an $m_0$-dependent \regvary\ \seq . A
typical example of such a process is a moving average process of order
$m_0\ge 1$ (MA$(m_0)$) given by
\beao
X_t= Z_t+\theta_1\,Z_{t-1}+\cdots+ \theta_{m_0}\,Z_{t-m_0}\,,\quad t\in\bbz\,,
\eeao
where $(Z_t)$ is an iid \regvary\ \seq\ with index $\alpha>0$.
Condition ${\bf RV}_\alpha$ is
straightforward since $(Z_t)$ is \regvary\ with limiting \ms s
concentrated on the axes. The \regvar\ of the  \fidi s of  $(X_t)$
is then an application of the continuous mapping theorem for \regvar ;
see Hult and Lindskog \cite{hult:lindskog:2005,hult:lindskog:2006}; cf. Hult et al.
\cite{hult:lindskog:mikosch:samorodnitsky:2005},
Jessen and Mikosch
\cite{jessen:mikosch:2006}.
\par
A related example is given by a \sv\ model
$X_t=\sigma_t\eta _t$, $t\in \bbz$, where $(\log \sigma_t)$
constitutes an MA$(m_0)$ process independent of
the iid \regvary\ \seq\ $(\eta_t)$ with index $\alpha$. If
$E\sigma^{\alpha+\epsilon}<\infty$ for some $\epsilon>0$ then $(X_t)$
is  \regvary\ with index $\alpha$; see e.g. Davis and Mikosch
\cite{davis:mikosch:2001,davis:mikosch:2009b}. By construction,
$(X_t)$ is $m_0$-dependent.
\par
For $m_0$-dependent \seq s the verification of the conditions of
Theorem~\ref{th:main} is simple.
\bpr\label{thm:m0depend}
Consider an $m_0$-dependent stationary \seq\ $(X_t)$ for some $m_0\ge
1$.
Assume that $(X_t)$  satisfies ${\bf  RV}_\alpha$ for some
$\alpha>0$ and $\E X=0$ if $\E |X|<\infty$.
Choose $b_n=n^{(1/\alpha)\vee 0.5+\delta}$ for any $\delta>0$.
Then
Theorem~\ref{th:main} holds with  $b_+= b_+(m_0+1)-b_+(m_0)$ in the
regions $\Lambda_n=(b_n,\infty)$.
\epr
\begin{proof}
Condition ${\bf AC}_\alpha$ is trivially satisfied for any choice of
constants $\delta_k\downarrow 0$ as $\kto$ and any sets
$\Lambda_n\subset (0,\infty)$ \st\ $n \P(|X|>b_n)\to 0$ as
$\nto$. Moreover, $b_+= b_+(m_0+1)-b_+(m_0)$ follows from
Bartkiewicz et al.
\cite{bartkiewicz:jakubowski:mikosch:wintenberger:2011}.
\par
It remains to prove that \eqref{eq:main} holds. In view of the
$m_0$-dependence of the \seq\ $(X_t)$ it is possible to split
the sum in \eqref{eq:main} into two sums of independent  subsums 
consisting of at most  $m_0$ summands. 
More precisely, with the convention that $X_j=0$ if $j>n$ we write
\beao
\overline S_n&=&\sum_{i=1}^n X_i\1_{\{|X_i|\le \delta_k \,x\}}\\
&=&\sum_{j=1,j\,\mbox{\tiny even}}^{[n/m_0]} \sum_{i=m_0 j+1}^{m_0(j+1)}X_i\1_{\{|X_i|\le \delta_k \,x\}}+\sum_{j=1,j\,\mbox{\tiny odd}}^{[n/m_0]} \sum_{i=m_0 j+1}^{m_0(j+1)}X_i\1_{\{|X_i|\le \delta_k \,x\}}\\
&=&\overline S_{n}'+\overline S_{n}''.
\eeao
Since 
\beao
\P(\overline S_n\ge \vep_k x)\le \P(\overline S_n'\ge \vep_k x/2)+\P(\overline S_n''\ge \vep_k x/2)\eeao 
we obtain an upper bound similar to \eqref{eq:main} 
but with sums of at most $[n/2 m_0]$ iid subsums.  
Therefore we may
assume without loss of generality that the $(X_t)$ in \eqref{eq:main} are iid. In view of
Remark~\ref{rem:0} and the conditions above we may assume
without loss of generality that the
summands in \eqref{eq:main} are mean corrected.
\par
For $\alpha\in (0,2)$, an application of Chebyshev's inequality and
Karamata's theorem
yield the estimate
\beao
\P\Big(\sum_{i=1}^n (X_i\1_{\{|X_i|\le \delta_k \,x\}}-\E X\1_{\{|X|\le \delta_k \,x\}})>\varepsilon_k
x\Big)
&\le& n\, (\vep_k x)^{-2} \E X^2 \1_{\{|X|\le \delta_k \,x\}}\\
&\sim & \vep_k^{-2} \delta_k^{2-\alpha} [n\,\P(|X|>x)]\,.
\eeao
Now choose e.g. $\delta_k=\ex^{-k}$
and $\vep_k=k^{-2}$. Then all assumptions on $(\vep_k)$ and
 $(\delta_k)$ in Theorem~\ref{th:main} are satisfied and
$\lim_{\kto }\vep_k^{-2} \delta_k^{2-\alpha} =0 $. Hence  \eqref{eq:main} is satisfied.
\par
In the case $\alpha>  2$, we use the Nagaev-Fuk inequality
(cf. Petrov \cite{petrov:1995}, p. 78, 2.6.5) for $p>\alpha$ and
Karamata's theorem as $\nto$, for $x \in\Lambda_n$:
\beao\lefteqn{
\P\Big(\sum_{i=1}^n (X_i\1_{\{|X_i|\le \delta_k \,x\}}-\E X\1_{\{|X|\le \delta_k \,x\}})>\varepsilon_k
x\Big)}\\&\le & c\,(\vep_k x)^{-p}
n\,\E |X|^p\1_{\{|X|\le x\delta_k\}}+ \ex^{-c (\vep_k x)^2/ n}\\
&\le &c\,\Big(\delta_k^{p-\alpha}\vep_k^{-p}+ \ex^{-c (\vep_k x)^2/ n}/
                                   [n\,\P(|X|>x)]\Big)[n\,\P(|X|>x)]\,.
\eeao
Choosing $\delta_k=\ex^{-k}$ and
$\vep_k=k^{-2}$, the requirements of
Theorem~\ref{th:main} are satisfied and
$\delta_k^{p-\alpha}\vep_k^{-p}$ becomes arbitrarily small for large
$k$. Moreover, $\sup_{x\in\Lambda_n}\ex^{-c (\vep_k x)^2/
  n}/[n\,\P(|X|>x)]\to 0$ by the choice of $(b_n)$. This proves
\eqref{eq:main} for $\alpha>2$.
\par
The boundary case $\alpha=2$ can be treated in a similar way by using
another version of the Nagaev-Fuk inequality; see Petrov \cite{petrov:1995},
p. 78, 2.6.4. We omit details.
\end{proof}

\subsection{Stochastic volatility model}\label{subsec:sv}
Consider a stationary \seq\ $(\sigma_t)$ of non-negative \rv s and
assume that $(Z_t)$ is an iid \seq\ which is independent  of
$(\sigma_t)$. The stationary \seq\
\beam\label{eq:sv}
X_t=\sigma_t\,Z_t, \quad t\in\bbz,
\eeam is
then called a {\em \sv\ model}. It is  a standard model in financial
\tsa ; see e.g. Andersen et
al. \cite{andersen:davis:kreiss:mikosch:2009}.

The main result of this section is a
large deviation principle for such models under
various assumptions.
\bth Consider a \sv\ model \eqref{eq:sv}  \st\ $Z$ is regularly varying
with index $\alpha>0$, $\E Z=0$ for $\alpha>1$ and
$\E\sigma_0^{2\alpha}<\infty$. Moreover, consider the following
additional conditions:
\begin{enumerate}
\item
$Z$ is symmetric.
\item
$\E\sigma_0^{p}<\infty$ for some $p>2\alpha$
and $(\sigma_t)$ is strongly mixing
with rate $(\alpha_j)$ \st\ $\alpha_j\le c j^{-a}$ for some $a>1$.
\end{enumerate}
The large deviation principle \eqref{eq:1} holds with
$b_+=\lim_{x\to\infty}{\P(Z>x)}/{\P(|Z|>x)}$ in the regions
$\Lambda_n=(b_n,\infty)$ under the following conditions:
\begin{itemize}
\item $0<\alpha<1$: $b_n=n^{\vep+1/\alpha}$ for any $\vep>0$\,.
\item
$1< \alpha<2$: Assume $(1)$ or $(2)$,
$b_n=n^{\vep+1/\alpha}$ for any $\vep>0$.
\item $\alpha>2$: Assume $(2)$ for some
  $a>\max(1,(\alpha-2)p/(2p-\alpha))$,
$b_n=\sqrt{n\log n} s_n$ for any \seq\ $(s_n)$
\st\ $s_n\to\infty$.
\end{itemize}
\ethe
\bre A Gaussian stationary process $(Y_t)$ is
strongly mixing under mild conditions; see
 Kolmogorov and Rozanov \cite{kolmogorov:rozanov:1960}.
Ibragimov \cite{ibragimov:1970}, Theorem~5,
gave necessary and sufficient conditions
for the relation $\alpha_n=O(n^{-a})$ for any choice of $a>0$.
The conditions are in terms of the spectral density of $(Y_t)$.
It is also known that a linear Gaussian process $Y_t=\sum_{j=0}^\infty
\psi_j\eta_{t-j}$, $t\in\bbz$, with $(\eta_t)$ iid standard normal and
exponentially decaying coefficients $(\psi_j)$ has an
exponentially decaying mixing rate $(\alpha_j)$; see
Pham and Tran \cite{pham:tran:1985}; cf. Doukhan \cite{doukhan:1994}.
For example,
if $(Y_t)$ is a causal Gaussian ARMA process
the latter condition is satisfied.
\par
Now assume $\log \sigma_t=Y_t$,
$t\in\bbz$, for a Gaussian stationary \seq\ $(Y_t)$.
This Gaussian model is chosen in the majority of the
literature on \sv\ models; see e.g. Andersen et al.
\cite{andersen:davis:kreiss:mikosch:2009}. Then $(\sigma_t)$
inherits strong mixing  from $(Y_t)$ with the same rate.
Of course, $\E\sigma^p<\infty$ for all $p>0$ and the large deviation principle holds for $\alpha_n=O(n^{-a})$ for any $a>1$.
\ere
\bre
If $(\sigma_t)$ is strongly mixing with rate $(\alpha_j)$,
the corresponding \sv\ model $(X_t)$
is strongly mixing with rate $(4\alpha_j)$; see e.g. Davis and Mikosch
\cite{davis:mikosch:2001}. 
\ere

\begin{proof} Condition ${\bf RV}_\alpha$ was verified for \sv\ models
under the condition $\E \sigma^{\alpha+\epsilon}<\infty$ for some $\epsilon>0$
in Davis and Mikosch \cite{davis:mikosch:2001}; see also
\cite{davis:mikosch:2009b}.
The limit \ms s of the \regvary\ \fidi s are concentrated on the axes
and therefore  $b_+=\lim_{x\to\infty}{\P(Z>x)}/{\P(|Z|>x)}$; see also
Bartkiewicz et al. \cite{bartkiewicz:jakubowski:mikosch:wintenberger:2011}.
\par
Next we verify condition ${\bf AC}_\alpha$. Fix any $\delta>0$. We have
\beao
p_{j}(\delta)=\P(|X_j|> x\delta,|X_0|> \delta x)\le \P( |Z_j Z_0|
\sigma_j\sigma_0)>(\delta x)^2)\,.
\eeao
The \rv\ $|Z_j Z_0|$ is \regvary\ with index
$\alpha$; see Embrechts and Veraverbeke \cite{embrechts:veraverbeke:1982}.
An application of Markov's and H\" older's inequalities yields for $\epsilon<2\alpha$,
\beao
p_{j}(\delta)\le   (\delta x)^{-2\alpha+\epsilon} (\E |Z|^{\alpha-\epsilon/2})^2
\,E|\sigma_j\sigma_0|^{\alpha-\epsilon/2}\le (\delta x)^{-2\alpha+\epsilon}
(\E |Z|^{\alpha-\epsilon/2})^2
\,E|\sigma|^{2\alpha-\epsilon}\,.
\eeao
We also have for any small $\epsilon>0$ and large $x$, $P(|X|>\delta
x)\ge (\delta x)^{-\alpha -\epsilon}$ in view of the \regvar\ of $X$. Therefore
\beao
 \sup_{x> b_n}\delta^{-\alpha}\sum_{j= k}^n\P(|X_j|> x\delta \mid |X_0|> x\delta)\le c\,n\,\delta^{-\alpha+2\epsilon}b_n^{-\alpha+2\epsilon}\,.
\eeao
The \rhs\ converges to zero if we choose $\alpha\le 2$,
$b_n=n^{\vep +1/\alpha}$ for any $\vep>0$ or $\alpha>2$,
$b_n=\sqrt{n\log n} s_n$, $s_n\to\infty$ and $\epsilon$ sufficiently small. The choice of
$\delta=\delta_k\to 0$ is arbitrary.
\par
Next we prove condition \eqref{eq:main}.\\
{\em The case $\alpha<1$.} Condition  \eqref{eq:main} is immediate
from Remark~\ref{rem:less1} for
$\delta_k=\ex^{-k}$ and $\vep_k=k^{-2}$.\\
The following decomposition will be useful in the case $\alpha>1$:
\beao
\lefteqn{\P\Big(\sum_{i=1}^n \sigma_i Z_i\1_{\{|\sigma_iZ_i|\le \delta_k
  \,x\}}>\vep_k x\Big)}\\&\le& \P\Big(\sum_{i=1}^n [\sigma_i Z_i\1_{\{|\sigma_iZ_i|\le \delta_k
  \,x\}}- \sigma_i \E (Z\1_{\{|\sigma_iZ|\le \delta_k
  \,x\}}\mid \sigma_i)]>(\vep_k/2) x\Big)\\&&+\P\Big(\sum_{i=1}^n \sigma_i \E (Z\1_{\{|\sigma_iZ|\le \delta_k
  \,x\}}\mid\sigma_i) >(\vep_k/2) x\Big)
=I_1+I_2\,.
\eeao
\ble\label{lem:help} Assume $\alpha>1$, and either $Z$ is symmetric or
$(\sigma_t)$ is strongly mixing with rate \fct\
$(\alpha_j)$ satisfying $\alpha_j\le c j^{-a}$ for some $c>0,a>1$ and
$E\sigma^p<\infty$ for some $p>2\alpha$.
Then
\beao
\lim_{\nto}\sup_{x>b_n}\dfrac{I_2}{n\,\P(|X|>x)}=0\,.
\eeao
\ele
\begin{proof}  In the case of symmetric $Z$, $I_2=0$.
Thus we deal with the case of mixing $(\sigma_t)$.
First observe that for any $y>0$,
\beao
I_2&\le& \P\Big(\sum_{i=1}^n \sigma_i \1_{\{\sigma_i\le y \}}
\E (Z\1_{\{|\sigma_iZ|\le \delta_k
  \,x\}}\mid\sigma_i) >(\vep_k/2) x\Big)+ n\,\P(\sigma>y)\\&=&I_{21}+I_{22}\,.
\eeao
Clearly, since $\E \sigma^{p}<\infty$ for some $p>2\alpha$, we can find
$y=y(x)=o(x)$, $y\to\infty$
as $\xto$ \st
\beao
\sup_{x>b_n} \dfrac{I_{22}}{n\,\P(|X|>x)}=\sup_{x>b_n}\dfrac{\P(\sigma>y)}{\P(|X|>x)}=o(1)\,.
 \eeao
Indeed, we can choose $y=x^{0.5-\gamma}$ for any $\gamma>0$ close to zero.
Write
\beao
\ov\sigma_i= \sigma_i \1_{\{\sigma_i\le y \}}
\E (Z\1_{\{|\sigma_iZ|\le \delta_k
  \,x\}}\mid\sigma_i)\,,\quad i=1,2,\ldots\,,
\eeao
and $\ov S_n=\sum_{i=1}^n \ov \sigma_i$.
The Markov inequality yields
\beao
\P (\ov S_n>\varepsilon_k x)&\le& (\varepsilon_k x)^{-2}
\E\ov  S_n^2\\
&=& (\varepsilon_k x)^{-2}\Big[ n \E \ov \sigma^2+2\sum_{j=1}^{n-1}(n-j) \E(\ov \sigma_0
\ov\sigma_j)\Big]
= I_3+ I_4\,.
\eeao
Then, since $\E Z=0$, by Karamata's theorem
\beao
\dfrac{I_3}{n\,\P(|X|>x)}\le
c\,\dfrac{x^{-2}[\E |Z|\1_{\{|Z|>\delta_k
    x/y\}}]^2}
{\P(|X|>x)}\le c\,\dfrac{y^{-2} [\P(|X> x/y)]^2}
{\P(|X|>x)}\,.
\eeao
The \rhs\ is negligible uniformly for $x>b_n$.
We also have
\beao
\dfrac{(n/x)^2 \,(\E \ov \sigma)^2}{n\P(|X|>x)}&=&
\dfrac{ n(\E (X_1 \1_{\{|X_1|>\delta_k
  x\,,\sigma_1\le y \}}))^2}{x^2 \P(|X|>x)}\\&\le&
\dfrac{n(E|X|\1_{\{|X|>\delta_k x\}})^2}{x^2\P(|X|>x) }\\
&\le& c\,n\,P(|X|>x)\le n\,P(|X|>b_n)\to 0 \,.
\eeao
Therefore we may assume without loss of generality that the
\rv s $\ov \sigma_j$ in $I_4$ are centered. Using a classical bound
for the covariance of a strongly mixing \seq , the fact that $\E
Z=0$ and Karamata's theorem,  for $r,q>0$ \st\ $r^{-1}+ 2q^{-1}=1$,
$1<r<a$,
\beao
|\cov(\ov \sigma_0,\ov\sigma_j)| &\le&
c\,\alpha_j^{1/r} \,(\E \ov \sigma^q)^{2/q}\\
&\le & c\,\alpha_j^{1/r} y^2[\E ( |Z|\1_{\{|Z|>\delta_k x/y\}})]^2\\
&\le & c\,\alpha_j^{1/r} x^{2}[\P(|Z|>x/y)]^2\,.
\eeao
Finally, we get the following bound
\beao
\sup_{x>b_n}\dfrac{I_4}{n\,\P(|X|>x)}\le
c\sum_{j=1}^\infty\alpha_j^{1/r}\,\sup_{x>b_n}\,
\dfrac{[\P(|Z|>x/y)]^2}{\P(|X|>x)}\,.
\eeao
The \rhs\
converges to zero. This proves the lemma.
\end{proof}\noindent
{\em The case $\alpha\in (1,2)$.} In view of Lemma~\ref{lem:help} it
remains to bound $I_1$.
Applying Chebyshev's inequality conditionally on $(\sigma_i)$ we
obtain
\beao
I_1&\le& (\vep_k x)^{-2} \E \Big[\sum_{i=1}^n \sigma_i^2 \var(Z\1_{\{|\sigma_iZ|\le \delta_k
  \,x\}}\mid \sigma_i)\Big]\\
&\le & (\vep_k x)^{-2} n\,\E(X^2\1_{\{|X|\le \delta_k
  \,x\}})\,.
\eeao
Now an application of Karamata's theorem and  \regvar\ of $X$ yield
\beao
\sup_{x>b_n}\dfrac{I_1}{n\,\P(|X|>x)}\le c\,\sup_{x>b_n} \dfrac{\delta_k^{2}}{\vep_k^2} \dfrac{\P
(|X|>x\delta_k)}{\P(|X|>x)}\sim c\,\dfrac{\delta_k^{2-\alpha}}{\vep_k^2}\,.
\eeao
Now choose $(\delta_k)$ and $(\vep_k)$ as in the case $\alpha<1$ to
conclude that
\beao
\lim_{\kto}\sup_{x>b_n}\dfrac{I_1}{n\,\P(|X|>x)}=0\,.
\eeao
The finishes the proof of \eqref{eq:main} in the case $\alpha\in (1,2)$.\\[1mm]
{\em The case $\alpha>2$.} We again have to study $I_1$.
Using the Nagaev-Fuk inequality (cf. Petrov \cite{petrov:1995}, p. 78, 2.6.5)
conditionally on $(\sigma_t)$, we obtain for $p>\alpha$,
\beao
\lefteqn{\P\Big(\sum_{i=1}^n
[\sigma_iZ_i \1_{\{|\sigma_i Z_i |\le \delta_k x\}}-\sigma_i\,\E(Z_i
\1_{\{|\sigma_i Z_i |\le \delta_k x\}}\mid
\sigma_i)]>(\vep_k/2)x\mid (\sigma_i)\Big)}\\
&\le& c\,(\vep _kx)^{-p}\sum_{i=1}^n\sigma_i^p\E(|Z_i|^p \1_{\{|\sigma_i Z_i |\le \delta_k x\}}\mid\sigma_i)+\ex^{-c(\vep_k x)^2/\sum_{i=1}^n\sigma^2_i}.\\
\eeao
The expectation of the first term is
of the \asy\ order $c\delta_k^{p-\alpha}/\vep_k^p$. The latter
relation converges to zero for $\delta_k=\ex^{-k}$ and $\vep_k=k^{-2}$.
Consider the expectation of the second term on the sets
$\{\sum_{i=1}^n\sigma^2_i > c(\vep_k x)^2/(2\alpha \log x)\}$
and its complement to obtain the bound
\beao
\E(\ex^{-c(\vep_k x)^2/\sum_{i=1}^n\sigma^2_i})\le
x^{-2\alpha}+\P\Big(\sum_{i=1}^n\sigma^2_i>c(\vep_k x)^2/(2\alpha\log x)\Big).
\eeao
The first term is negligible with respect to $n\P(|X|>x)$.
For the second one, note that  $x^2/(n\log x)\ge c\, b_n^2/(n\log b_n)\to \infty$. Therefore we
may assume without loss of generality that the $\sigma_i^2$'s are mean
corrected. Now use Rio \cite{rio:2000}, p. 87, (6.19a),
under the mixing condition
$\alpha_j\le cj^{-a}$ to obtain for any
$r\ge 1$:
\beao\lefteqn{\P\Big(\sum_{i=1}^n
(\sigma^2_i-\E\sigma^2)>c(\vep_k x)^2/(2\alpha\log x)\Big)}\\
&\le& c\,n^{r/2}(\log x)^{r}x^{-2r}+c\,n(\log(x)/x^2)^{(a+1)p/(a+p)}.
\eeao
The first term is negligible with respect to $n\,P(|X>x)$ for $r$
sufficiently large. The second term is negligible as well
if $2(a+1)p/(a+p)>\alpha$. The latter condition is satisfied by assumption.
\end{proof}

\subsection{Regularly varying functions of Markov chains} \label{subsec:mc}
In this section we assume that { $X_t=h(\Phi_t)$, $t\in\bbz$,} is 
a { measurable real-valued function} of a stationary \MC\  $(\Phi_t)$ which possesses an atom $A$ in some general space: 
The context is classical; see Nummelin
\cite{nummelin:1984}
and Meyn and Tweedie
\cite{meyn:tweedie:1993} which will serve as our main references, and $(\Phi_t)$ can be seen as the enlargement of a Harris recurrent \MC. In
Section~\ref{subsec:sre} we will look at the example of a
solution to a \sre\ which constitutes such a \MC . We assume that the function $h$ is such that $(X_t)$ is  \regvary\ with index $\alpha>0$. Notice in particular that { $h$ is not the null \fct .}
\par
Throughout we will also assume the following
{\em polynomial drift condition} for  $p>0$ which is inspired by Samur
\cite{samur:2004}
who used a more general condition.
\begin{itemize}
\item
 ${\bf DC}_p$:
There exist constants $\beta\in  (0,1)$, $b>0$  such that for any $y$,
\beao
\E( |h(\Phi_1)|^p\mid \Phi_0=y)\le \beta \,|h(y)|^p+b\,\1_A(y).
\eeao
\end{itemize}
In this condition, we suppress the dependence of $\beta,b,A$ on the
value $p$.
Note that ${\bf DC}_p$ implies geometric
ergodicity of $(\Phi_t)$; see Meyn and Tweedie \cite{meyn:tweedie:1993},
p. 371.
In what follows,
we write  $\tau_A$ for the first time the
chain visits the set $A$, $\P_A$
denotes the \pro y \ms\ of the \MC\ conditional on $\{\Phi_0\in A\}$ and
$\E_A$ is the corresponding expectation. We will also write $\P_x$ and
$\E_x$ if $ \{\Phi_0=x\}$.
\par
Here is the main result of this section.
\bth\label{thm:mc}
Assume that $(\Phi_t)$ is a stationary \MC\ possessing an atom $A$ and that $h$ is a function such that { $X_t=h(\Phi_t)$, $t\in\bbz$, satisfies} the conditions $(1)-(3)$ of
Theorem~\ref{th:main}  for the regions $\Lambda_n=(b_n,c_n)$ specified below.
Also assume $\E X=0$ if $\E |X|<\infty$ and  ${\bf DC}_{p}$ for all $p<\alpha$.
Then the precise \ld\ principle \eqref{eq:1} holds under the following
conditions:
\begin{itemize}
\item
$0< \alpha<1$: $\Lambda_n=(b_n,\infty)$ for any \seq\ $(b_n)$
satisfying $n\P(|X|>b_n)\to 0$.
\item
$1< \alpha$ and $\alpha\neq2$: $\Lambda=(b_n,c_n)$ for any  \seq\ $(b_n)$
satisfying $b_n=n^{1/\alpha\vee 0.5+\delta}$ for  any $\delta>0$,
and $(c_n)$ \st\ $c_n>b_n$ { and }
\beam\label{eq:tau}
\P(\tau_A>n)= o(n\, \P(|X|>c_n))\,.
\eeam
\end{itemize}
\ethe
\begin{proof}
We will apply Theorem~\ref{th:main}.
Since we assumed conditions
(1)--(3) of this result it remains to verify \eqref{eq:main}.\\[1mm]
{\em The case $0<\alpha<1$.} The proof follows from Remark~\ref{rem:less1}.\\
{\em The case $\alpha>1$ and $\alpha\notin \N$.} This case is more involved. We will prove
it in  a similar way as in the iid or $m_0$-dependent cases, by using
moment and exponential inequalities tailored for regenerative split \MC s. Without loss of generality we will only consider the strongly
aperiodic case. 

Notice that ${\bf DC}_p$ is satisfied for $p=[\alpha]$. For all integers $p<[\alpha]$, applying Jensen's inequality, we obtain 
\beam\label{eq:jensen}
\E(|X_1|^p\mid \Phi_0=y)&\le& \Big(\E(|X_1|^{[\alpha]}\mid \Phi_0=y)\Big)^
{p/[\alpha]}\nonumber\\&\le& \Big(\beta \,|h(y)|^{[\alpha]}+b\,\1_A(y)\Big)^{p/[\alpha]}\le \beta^{p/[\alpha]}\,|h(y)|^p+b^{p/[\alpha]}\,\1_A(y).
\eeam
Thus $b>0$, $\beta\in (0,1)$ 
and $A$ in ${\bf DC}_p$ can be chosen
the same as in  ${\bf DC}_{[\alpha]}$. 

Let $(\tau_A(j))_{j\ge 1}$ be the \seq\ of visiting times of the \MC\
to the set $A$, i.e. $\tau_A(1)=\tau_A$ and
$\tau_A(j+1)=\min\{k>\tau_A(j): \Phi_k\in A\}$. Notice that the \seq\
$(\tau_A(j+1)-\tau_A(j))_{j\ge 1}$ constitutes an iid \seq\ and
$N_A(t)=\#\{j\ge 1: \tau_A(j)\le t\}$, $t\ge 0$, is a renewal process.
The following inequality holds for any integrable \fct\ $f$ on $\bbr$:
\beao
\lefteqn{\P\Big(\sum_{i=1}^n f(X_i)> \vep_k x\Big)}\\
&=& \P\Big(\sum_{i=1}^n f(X_i)> \vep_k x\,,N_A(n)=0\Big)
+\P\Big(\sum_{i=1}^n f(X_i)> \vep_k x\,,N_A(n)=1\Big)\\
&&+\P\Big(\sum_{i=1}^n f(X_i)> \vep_k x\,,N_A(n)\ge 2\Big)\\
&\le &
\P( \tau_A>n)+2\P\Big(\sum_{j=1}^{\tau_A}
f(X_j)>\vep_kx/3, \tau_A\le n \Big)\\
&&+\P\Big(\sum_{j=1}^{N_A(n)-1}
\sum_{t=\tau_A(j)+1}^{\tau_A(j+1)}f(X_j)>\vep_kx/3\Big)+2\P\Big(
\sum_{i=\tau_A(N_A(n))+1}^{n}f(X_i)>\vep_k x/3\Big)\\
&=& I_1+I_2+I_3+I_4\,.
\eeao
We mentioned in Remark~\ref{rem:0} that we may assume without loss of
generality that the \rv s $\ov X_i$, $i=1,2,\ldots,$ are mean
corrected. Now we choose $f(X_i)=\ov X_i -\E \ov X_i $ where
\beao
\ov X_i=X_i\1_{\{|X_i|\le \delta_k x\}}\,,\quad i=1,2,\ldots,\quad x>0\,.
\eeao
{\em Bounds for $I_1,I_2,I_4$.} For $I_4$, we
use the Markov inequality of order  $k_0=[\alpha]+1$ and the   stationarity of
$(X_i)$ 
\beao
I_4&\le& c \,(x\vep_k)^{-k_0}\Big[
\E\Big| \sum_{i=\tau_A(N_A(n))+1}^{n}\ov X_i\Big|^{k_0}\Big] +\E
\tau_A^{k_0}\; [\E |X|\1_{\{|X|>\delta_k x\}}]^{k_0}\Big] \\
&\le& c\,x^{-k_0}\Big[
\E\Big( \sum_{i=\tau_A(N_A(n))+1}^{n}\,|\overline X_i|\Big)^{k_0}+ [x\,\P (|X|>x)]^{k_0}
\Big]\\
&\le& c\,x^{-k_0}\Big[
\E_A\Big( \sum_{i=1}^{ \tau_A }|\overline X_i|\Big)^{k_0}+[x\,\P (|X|>x)]^{k_0}\Big] \, .
\eeao
Since for $\alpha>1$, $k_0\ge 2$, we  use
 Proposition \ref{pr:mineq} given below 
to show that $I_4$ is negligible
with respect to $n\P(|X|>x)$.
As to  $I_2$, we again use the Markov inequality: 
\beao
I_2&\le& c\, (x\vep_k)^{-k_0}\Big[
\E\Big|1_{\{\tau_A\le n\}} \sum_{i=1}^{\tau_A }\ov X_i\Big|^{k_0}
+\E \tau_A^{k_0}\; [\E |X|\1_{\{|X|>\delta_k x\}}]^{k_0}\Big]\\
&\le& c\,x^{-k_0}\Big[
\E\Big(1_{\{\tau_A\le n\}} \sum_{i=1}^{\tau_A }|\overline
X_i|\Big)^{k_0}
+ [x\,\P (|X|>x)]^{k_0}
\Big]
\,.
\eeao
We
iteratively apply Lemma~\ref{lem:2} given below to the first term in the \rhs\ to obtain an estimate
of $I_2$ proportional to
\begin{equation}\label{eq:pt}
x^{-k_0}\E\Big(1_{\{\tau_A\le n\}}  \sum_{i=1}^{  \tau_A }|\overline X_i|^{k_0}\Big)=x^{-k_0}\E\Big(  \sum_{i=1}^{   n }|\overline X_i|^{k_0}\1_{\{\tau_A\ge i\}}\Big).
\end{equation}
An application of Pitman's identity \cite{pitman:1977} yields
\beao
\E\Big(  \sum_{i=1}^{   n }|\overline X_i|^{k_0}\1_{\{\tau_A\ge i\}}\Big)&=& \P(\Phi_0\in A)\,\E_A\Big(  \sum_{k=0}^{\tau_A-1}\sum_{i=1}^{   n }|\overline X_{k+i}|^{k_0}\1_{\{\tau_A\ge k+i\}}\Big)\\
&\le& n  \P(\Phi_0\in A)\,\E_A\Big(\sum_{i=1}^{\tau_A}|\overline X_{i}|^{k_0}\Big).
\eeao
From a Wald-type identity, $I_2\le cn(x\vep_k)^{-k_0}\E|\overline X|^{k_0}$.
Hence  $I_2$ is negligible with respect to $n\P(|X|>x)$ by an application of Karamata's theorem.
\par
Finally, $I_1$ is negligible   with respect to $n\P(|X|>x)$
because we assume that  $\P(\tau_A>n)=o(n\P(|X|>c_n))$.\\[1mm]
{\em Bounds for $I_3$.}
The following moment inequality is the key to the bound of
$I_3$:
\bpr\label{pr:mineq}
Assume that { $(X_t)=(h(\Phi_t))$ for a real-valued measurable \fct\ $h$ 
and a Markov chain $(\Phi_t)$} satisfying the drift condition
${\bf DC}_{k_0-1}$ for some integer $k_0\ge
2$. Then for $x>0$ and some constant $c>0$,
\beam\label{eq:11}
\E_{A}\Big(\sum_{j=1}^{\tau_A} |\overline
X_j|\Big)^{k_0}\le c\,\E |\overline X|^{k_0}\,.
\eeam
\epr
\begin{proof} We can expand the \lhs\ of \eqref{eq:11} as follows
\beam\label{eq:tetra}
\E_A\Big(\sum_{j=1}^{\tau_A} |\overline
X_j|\Big)^{k_0}=\sum_{k=1}^{k_0}\sum_{\sum_{i=1}^ks_i=k_0,s_i\ge
  1,i=1,\ldots,k}\E_A\Big(\sum_{j_1=1}^{\tau_A}
\sum_{j_2=j_1+1}^{\tau_A}\cdots\sum_{j_k=j_{k-1}+1}^{\tau_A} |\overline X_{j_i}|^{s_i}\Big)\,.
\eeam
We will estimate the moments on the \rhs\ by employing
Lemma~\ref{lem:2} below. For the cases $k_0=2,3$ such a result was
proved by  Samur \cite{samur:2004} and we use the idea of the proof in
\cite{samur:2004} for our generalization.
Before we formulate the basic
moment estimate we need some notation: According to
the proof of Theorem 14.2.3 of Meyn and Tweedie \cite{meyn:tweedie:1993}, there exists
a constant $c(A)>0$ \st
$$
\E_{\Phi_0}\Big(\sum_{k=1}^{\tau_A} 1_A(X_k)\Big)\le c(A)\quad \as
$$
\ble\label{lem:2}
Assume ${\bf DC}_{p}$ and let $f,g$ be non-negative measurable \fct s
on $\bbr$ \st\ $f(x)\le |y|^p$ and $g(y)=0$ for $|y|>\delta_k x$.
Then for any $\ell\ge 1$, $n\in\N\cup\{\infty\}$
\beqq\label{eq:cmineq}
\E\Big(1_{\{\tau_A\le n\}}\sum_{j=\ell}^{\tau_A }g(X_j)\sum_{i=j+1}^{\tau_A }f(X_i)\mid
\bbf_{\ell}\Big)\le \E \Big(1_{\{\tau_A\le n\}}\sum_{j=\ell}^{\tau_A }g(X_j)[C \,|\overline
X_j|^p+  b\, c(A)]\mid \bbf_{\ell}\Big)\,,
\eeqq
where $\bbf_\ell=\sigma((\Phi_t)_{t\le \ell})$.
\ele
\begin{proof}
As mentioned in Samur \cite{samur:2004},
$\{\tau_A\ge j\}\in  \bbf_j$ for all $j$. Therefore
\begin{eqnarray*}
\E \Big(1_{\{\tau_A\le n\}}\sum_{j=\ell}^{\tau_A }g(X_j)\sum_{i=j+1}^{\tau_A }f(X_i)\mid
\bbf_{\ell}\Big)&=& \sum_{j=\ell}^n\E \Big(\1_{\{\tau_A\ge
  j\}}g(X_j)\sum_{i=j+1}^{\tau_A}f(X_i)\mid \bbf_{\ell}\Big)\\
&=& \sum_{j=\ell}^n\E\Big(\1_{\{\tau_A\ge
  j\}}g(X_j)\E\Big(\sum_{i=j+1}^{\tau_A}f(X_i)\mid \bbf_{j}\Big)\mid \bbf_{\ell}\Big)\\
&\le& \sum_{j=\ell}^n\E\Big(\1_{\{\tau_A\ge
  j\}}g(X_j)\;\E_{\Phi_j}\Big(\sum_{i=1}^{\tau_A}f(X_i)\Big)\mid \bbf_{\ell}\Big)\,.
\end{eqnarray*}
In the last inequality we used the stationarity  of $(\Phi_t)$ and
the strong Markov property. From Theorem 14.2.3 of Meyn and Tweedie
\cite{meyn:tweedie:1993} we obtain
$$
\E_{\Phi_j}\Big(\sum_{i=1}^{\tau_A}f(X_i)\Big)\le  C\,|X_j|^p+  b\,c(A).
$$
Since $g$ vanishes for $|y|\ge x$ the result for the truncated \rv s
$\ov X_j$ follows. This finishes the proof of Lemma~\ref{lem:2}.
\end{proof}
\noindent
By \eqref{eq:jensen} for $1\le p\le k_0-1$, 
${\bf DC}_p$ is satisfied for the same choice of $(b,A)$.  We can iteratively apply Lemma~\ref{lem:2} to
the expectations of the tetrahedral sums on the \rhs\ of
\eqref{eq:tetra}, starting with the tetrahedron with the largest index.  In the last step of the iteration we are left with a
sum of the type
\beao
\E_A\Big(\sum_{i=1}^{\tau_A}|\ov X_i|^{k_0})\Big)=\E |\ov X|^{k_0}\,\E_A(\tau_A)\,,
\eeao
where we used Wald's identity for any bounded $f$ on the \rhs . Thus, each of the summands
on the \rhs\ of \eqref{eq:tetra} can be bounded by the expression
$$
\E_A (\tau_A) \,\E|\overline X|^{k_0}\, \sum_{j=0}^k C^{k-j}(b \,c(A))^j
$$
and so the desired result follows.
\end{proof}
\noindent
{\em Bounds for $I_3$ in the case $1<\alpha<2.$}
 By  Markov's
inequality of order 2,
$$
\P\Big(\sum_{j=1}^{N_A(n)-1}
\sum_{t=\tau_A(j)+1}^{\tau_A(j+1)}f(X_j)>\vep_kx/3\Big)\le c(\vep x)^{-2}\E\Big(\sum_{j=1}^{N_A(n)-1}
\sum_{t=\tau_A(j)+1}^{\tau_A(j+1)}f(X_j)\Big)^2.
$$
From the regeneration scheme, we know that the  cycles
$(\sum_{t=\tau_A(j)+1}^{\tau_A(j+1)}f(X_j))$ are independent. Thus we
can expand the expectation
term and bound it by
$
n\E_A[S_A(f)^2].
$
The desired result follows by an application of
Proposition \ref{pr:mineq} with $k_0=2$
and Karamata's Theorem.\\[2mm]
{\em Bounds for $I_3$ in the case $\alpha>2$ and $\alpha\notin \N$.} The following inequality of Bertail and
Cl\'emencon \cite{bertail:clemencon:2009} is the key to the bound of
$I_3$ for $\alpha>2$.
It will be convenient to write $S_A(f)=\sum_{i=1}^{\tau_A} f(X_i)$.
\ble
Assume that $\sigma_A^2=\E_A \tau_A^2<\infty$ and
$\sigma_f^2=\E _A [(S_A(f))^2]<\infty$.
 Then for any $x$, sufficiently large $n$, $ \bfM=(M_1,M_2)\in (0,\infty)^2$ with Euclidean
norm $\|\bfM\|$,
\beam\label{eq:exp}
I_3&\le & c_0 \,\|\bfM\|^2\,
\exp\Big\{ - \dfrac{n(1+|\tilde \rho|)\tilde \sigma^2}
{2\|\bfM\|^2}
  H\Big(
\dfrac{\sqrt{2} \|\bfM\| \vep_k x} {n(1+|\tilde \rho|\tilde
    \sigma\tilde \sigma_f}\Big)\Big\}\\
&&+ (n-1) \P_A(|S_A(f)|>M_1)+ (n-1) \P_A(\tau_A>M_2)\,,
\eeam
where $H$ is the Bennett function $H(x)=(1+x)\ln(1+x)-x$,
$\tilde \sigma_f^2 =\var_A( S_A(f)\1 _{\{|S_A(f)|\le M_1\}})$,
$\tilde \sigma_A^2 =\var_A( \tau_A\1 _{\{|\tau_A|\le M_2\}})$,
$\tilde \rho=(\tilde\sigma_A\tilde \sigma_f)^{-1}\cov_A( S_A(f)\1
_{\{|S_A(f)|\le M_1\}}, \tau_A\1 _{\{|\tau_A|\le M_2\}})$, $\tilde
 \sigma^2 =\tilde \sigma_f^2\tilde\sigma_A^2/(\tilde \sigma_f^2+\tilde\sigma_A^2)$,
and some $c_0>0$.
\ele
Bertail and Cl\'emencon \cite{bertail:clemencon:2009} also assume that $\E_A S_A(f)=0$. This
condition is always satisfied in our situation since $\E f(X)=0$;
see Meyn and
Tweedie \cite{meyn:tweedie:1993}, (17.23) in Theorem 17.3.1. Under our
conditions, $\sigma_A^2$ is finite for any $\alpha$ and
$\sigma_f^2$ is finite for
$\alpha>2$; see Proposition~\ref{pr:mineq}.
One even has the stronger property: there exists a constant $\kappa>0$ \st
\beam\label{eq:tweedie}
\sup_{x\in A}\E_x\ex^{\kappa \tau_A}<\infty\,,
\eeam
see  Meyn and
Tweedie \cite{meyn:tweedie:1993},  (15.2) in Theorem~15.0.1.
We will choose
$M_1=M_2=\gamma_k\,x$ for some constants $\gamma_k>0$.
A careful study of the proof in \cite{bertail:clemencon:2009} shows
that $\tilde \rho,\tilde\sigma,\tilde \sigma_f$ are bounded for $\alpha>2$.
Then the exponential inequality \eqref{eq:exp} turns into
\beao
I_3
&\le& c \,(x\gamma_k)^2\ex^{-cn/(x\gamma_k)^2H(c x^2\gamma_k\vep_k/n)}
+n\,\P_A(\tau_A>
x\gamma_k)+n\,\P_A(|S_A(f)|>x\gamma_k)\nonumber\\
&=&I_{31}+I_{32}+I_{33}\,,
\eeao
for suitable constants $c>0$.
Choose
$\gamma_k=o(\varepsilon_k)$. Then for $k$ large,
uniformly for $x\ge b_n$ \st\ $b_n/n^{\delta+0.5}\to\infty$ for some
$\delta>0$,
\beao
\dfrac{I_{31}}{n\,\P(|X|>x)}\le
c\,\dfrac{x^{2(1-c\varepsilon_k/\gamma_k)}n^{c\varepsilon_k/\gamma_k}}{n\,\P(|X|>x)}=o(1)\,,\quad
\nto\,.
\eeao
As for $I_{32}$, it follows from \eqref{eq:tweedie}
and Markov's inequality that
\beao
\dfrac{I_{32}}{n\,\P(|X|>x)}\le
c\,\dfrac{\ex^{-\kappa x\gamma_k}}{\P(|X|>x)}
=o(1)\,,
\eeao
uniformly for $x\ge b_n$. Finally, Markov's inequality, an application of
Proposition \ref{pr:mineq} to $I_{33}$ with $k_0=[\alpha]+1$ and
Karamata's theorem
yield
\beao
\P_A(|S_A(f)|>x\gamma_k)\le  (x\gamma)^{-k_0}
\E| S_A(f)|^{k_0}\le c\,
 (x\gamma_k)^{-k_0}
\E| \overline X|^{k_0}\sim c\,\delta_k^{k_0-\alpha}\gamma_k^{-k_0} \P(|X|>x)\,.
\eeao
Choose $\delta_k=o(\gamma_k^{k_0/(k_0-\alpha)})$ as $\kto$.
This is always possible because we may choose $\vep_k=k^{-2}$,
$\delta_k=\ex^{-k}$ and $\gamma_k=k^{-3}$ throughout the proof.
Then we obtain
$$
\lim_{\kto}\sup_{x\ge b_n}\dfrac{I_{33}}{n\,\P(|X|>x)}\le
\lim_{\kto} c\,\delta_k^{k_0-\alpha}\gamma_k^{-k_0}=0\,.
$$
Thus we proved for $\alpha>2$
\beao
\limsup_{\kto}\lim_{\kto}\sup_{x\ge b_n}\dfrac{I_{3}}{n\,\P(|X|>x)}=0\,.
\eeao

{\em The case $\alpha>2$ and $\alpha\in \N$.} In this case, let us fix $\alpha/(\alpha+1)<\beta<1$ and consider the process $(|X_t|^\beta=|h(\Phi_t)|^\beta)$. It satisfies {\bf DC$_\alpha$} and concavity of $x\to x^\beta$ as $\beta<1$ implies that $$\E_A\Big(\sum_{i=1}^{\tau_A}|\overline X_i|\Big)^{\beta k_0}\le \E_A\Big(\sum_{i=1}^{\tau_A}|\overline X_i|^\beta\Big)^{k_0}.$$
We apply Proposition \ref{pr:mineq} to $(|X_t|^\beta)$ with $k_0=\alpha+1$ and we obtain  $\E_A|\overline S_A|^{\beta k_0}\le \E|\overline X_1|^{\beta k_0}$. Noticing that $\beta k_0>\alpha$, the use of Karamata's theorem as above yields that $\E|\overline X_1|^{\beta k_0}$ is negligible with respect to $n\P(|X|>x)$. Now we can follow the lines of the proof in the case of non-integer $\alpha$.
\end{proof}

In what follows, we will use the notation of
Theorem~\ref{thm:mc}
and its proof.
Our next goal is to give an intuitive interpretation
of the \ld\ principle of Theorem~\ref{thm:mc}: we want to show that
the \ld\ \pro y $\P(S_n>x)$ is essentially
determined by $\P(\max_{i=1,\ldots,N_A(n)} S_{A,i}>x)$, where
\beao
S_{A,i}=\sum_{t=\tau_A(i)+1}^{\tau_A(i+1)} X_t\,, \quad i\in\bbz\,,
\eeao
and
$(N_A(t))_{t\ge 0}$ is the renewal process generated from the iid \seq\ $(\tau_A(j+1))-\tau_A(j))$. 
The \seq\ $(S_{A,i})$ constitutes an
iid \seq . We write $\tau_A=\tau_A(1)$, $S_A=\sum_{i=1}^{\tau_A}X_i$ and $\lambda=(\E \tau_A)^{-1}$.
\bth\label{thm:int}
Assume that the  conditions of Theorem~\ref{thm:mc} hold, $\alpha>1$, $\alpha\neq 2$ and $b_+>0$.
Then $\P_A(S_A>x)\sim \E(\tau_A)b_+\P(|X|>x)$ and the precise \ld\ principle for the function of \MC\ $(X_t)$
can be written in the form
\beao
\sup_{x\in \Lambda_n}\Big|\dfrac{\P(S_n>x)}{n\,\P_A(S_A>x)}-
(\E \tau_A)^{-1}\Big|\to 0\,,
\eeao
where $\Lambda_n=(b_n,c_n)$ is chosen as in
Theorem~\ref{thm:mc}.
\ethe
\begin{proof}
Using the disjoint partition $\{N_A(n)=0\},\{N_A(n)=1\},\{N_A(n)\ge 2\}$, 
we obtain 
\beao
\P(S_n>x)&=& \P\Big( \sum_{i=1}^{n} X_i>x,\tau_A>n\Big)+ 
\P\Big( \sum_{i=1}^{\tau_A(1)} X_i+  \sum_{i=\tau_A(1)+1}^{n} X_i>x,\tau_A(2)>n\ge \tau_A(1)\Big)\nonumber\\&&+ \P(S_n>x, N_A(n)\ge 2)\nonumber.
\eeao
Using the definitions of $(\tau_A(i))$ and $N_A(n)$, we obtain for 
small $\vep\in (0,1)$
\beao
\P(S_n>x)&\le&\P( \tau_A>n)+2 \P(S_A>x\vep/2,\tau_A\le n )+\P \Big( \sum_{i=1}^{N_A(n)-1} S_{A,i}>x (1-\vep)\Big)\\
&&+2\P\Big( \sum_{i=\tau_A(N_A(n))+1}^{n} X_i> x \vep/2\Big)=J_1+J_2+J_3+J_4.
\eeao
and 
\beao
\P(S_n>x)&\ge& \P\Big(S_A+ \sum_{i=1}^{N_A(n)-1}S_{A,i}+ \sum_{t=\tau_A(N_A(n))+1}^n X_t>x\,,N_A(n)\ge 2\Big)\\
&\ge &\P\Big(\sum_{i=1}^{N_A(n)-1}S_{A,i}\ge (1+\vep) x,|S_A| \le \vep x/2,\Big|
\sum_{t=\tau_A(N_A(n))+1}^n X_t\Big|\le \vep x/2\,,N_A(n)\ge 2\Big)\\
&\ge &\P\Big(\sum_{i=1}^{N_A(n)-1}S_{A,i}\ge (1+\vep) x\Big)
-\P(|S_A|>\vep x/2)\\
&&-
\P\Big(\Big|\sum_{t=\tau_A(N_A(n))+1}^n X_t\Big|>\vep x/2\Big)-\P(N_A(n)\le 2)\\
&=&J_5-J_6-J_7-J_8\,.
\eeao
\ble\label{lem:1} Under the conditions of the theorem, for any small $\vep>0$, uniformly for $x\in \Lambda_n$, 
\beao
\dfrac{\P \Big( \sum_{i=1}^{N_A(n)-1} S_{A,i}>x (1+\vep)\Big)}{n\P(|X|>x)} +o(1)
\le \dfrac{\P(S_n>x)}{n\P(|X|>x)}\le 
\dfrac{\P \Big( \sum_{i=1}^{N_A(n)-1} S_{A,i}>x (1-\vep)\Big)}{n\P(|X|>x)} +o(1)\,.
\eeao
\ele
\begin{proof}
By assumption, the \pro y $J_1\le \P(\tau_A>n)$ is negligible with respect to $n\P(|X|>x)$
on $\Lambda_n$.

By standard computations and using the same notation as in the proof of Theorem \ref{thm:mc} we have
$$
J_4/2\le \P\Big( \sum_{i=\tau_A(N_A(n))+1}^{n} \overline X_i> x \vep/2\Big)+\P\Big(\cup_{i=\tau_A(N_A(n))+1}^{n}\{|X_i|>x\delta\}\Big).
$$
The second term is estimated by 
$$\E(\sum_{i=\tau_A(N_A(n))+1}^{n}1_{\{|X_i|>x\delta\}})\le \E_A(\sum_{i=1}^{\tau_A}1_{\{|X_i|>x\delta\}}) =\E(\tau_A)\,\P(|X|>x\delta).$$ 
The first term can be shown to be negligible with respect to $n\P(|X|>x)$ as in the proof of Theorem \ref{thm:mc}. So $J_4=o(n\P(|X|>x))$.

{ The term $J_2$ can be  treated in the same way  as  $I_2$ in the proof of Theorem~\ref{th:main}.} An application of Markov's 
inequality yields an estimate of the form $cx^{-k_0}\Big[\E\Big(1_{\{\tau_A\le n\}}  \sum_{i=1}^{  \tau_A }|\overline X_i|^{k_0}\Big)+ [{ n}\,\P (|X|>x)]^{k_0}
\Big]$. Using \eqref{eq:pt}, Pitman's and Wald-type identities we obtain $J_2\le cn(x\vep)^{{ -k_0}} \E|\overline X|^{k_0}$. Hence  $J_2$ is negligible with respect to $n\P(|X|>x)$ by an application of Karamata's theorem.

Collecting the 
bounds above, the upper bound in the lemma is proved.
\par
As regards the lower bound, $J_6$ and $J_7$ are of the order $o(n\P(|X|>x))$ in 
view of the bounds for $J_2$ and $J_4$ in the proof above, respectively. Moreover,
\beao
J_8= \P(N_A(n)\le 2)&\le& \P(\tau_A>n)+\P(\tau_A(2)>n)
\le 3 \,\P(\tau_A>n/2)\,,
\eeao
and the latter \pro y is negligible with respect to $n\P(|X|>x)$ as for $J_1$
above.
\end{proof} 

Denote $\tilde \Lambda_n=(b_n,\ex^{s_n})\cap \Lambda_n$
for some $(s_n)$ \st\ $s_n/n\to0$.
 \ble\label{lem:3} Under the conditions of the theorem, for any small $\xi,\vep>0$, uniformly for $x\in \tilde\Lambda_n$, 
\beam
\label{eq:lb}\lefteqn{\dfrac{\la (1-\vep)\P  (
S_{A}>x(1+\xi)(1+\vep) )}{\,\P(|X|>x)}+o\Big(\dfrac{\P_A( S_{A}>x )}{\P(|X|>x)}\Big)+o(1)}\\
&&\nonumber
\le \dfrac{\P \Big( \sum_{i=1}^{N_A(n)-1} S_{A,i}>x \Big)}{n\P(|X|>x)}\\
&&\label{eq:ub}\le 
\dfrac{\lambda\P_A( S_{A}>x (1-\xi))}{\P(|X|>x)} +o\Big(\dfrac{\P_A( S_{A}>x )}{\P(|X|>x)}\Big)+o(1)
\,.
\eeam
\ele
\begin{proof}
We have for $\delta>0$,
\beao
\P \Big( \sum_{i=1}^{N_A(n)-1} S_{A,i}>x\Big)
&=&\P \Big( \sum_{i=1}^{N_A(n)-1} S_{A,i}>x\,,|N_A(n)-1-n\,\la|>\delta n\Big)\\
&&+\P \Big( \sum_{i=1}^{N_A(n)-1} S_{A,i}>x\,,|N_A(n)-1-n\,\la|\le\delta
n\Big)\\
&=&K_{1}+K_{2}\,.
\eeao
In view of \eqref{eq:tweedie}, $\tau_A$ has exponential moment and therefore
one can apply standard \ld\ theory (e.g. Cram\'er's theorem; see Dembo
and Zeitouni \cite{dembo:zeitouni:2010}) to obtain
\beao
K_{1}\le \P(|N_A(n)-1-n\,\la|>\delta n) \le \ex^{-\gamma n}\,,
\eeao
for some $\gamma=\gamma(\delta)>0$. In view of the definition of $\tilde \Lambda_n$, $K_1=o(n\,\P(|X|>x))$ on $\tilde\Lambda_n$.
We also have
\beao
 \P \Big(
\sum_{i=1}^{n\la}  S_{A,i}-
\max_{|m-n \la |\le \delta n}\Big|\sum_{i=m}^{n\la} S_{A,i}\Big|>x\Big)
\le K_{2}\le \P \Big(
\sum_{i=1}^{n\la}  S_{A,i}+
\max_{|m-n \la |\le \delta n}\Big|\sum_{i=m}^{n\la} S_{A,i}\Big|>x\Big)\,.
\eeao
Here we define $\sum_{i=m}^b$ for any real value $b\ge m$, $m\in \bbn$, as
$\sum_{i=m}^{[b]}$ and the sums  $\sum_{i=b}^m$ are defined
accordingly.
Notice that  $b_n^{-1}\sum_{i=1}^{n\la}  S_{A,i}\stp 0$ from the fact that $n^{-1}N_A(n)\stas \la$.
Then, for any $\xi\in (0,1)$, a maximal inequality of
L\'evy-Ottaviani-Skorokhod type
for sums of iid \rv s (e.g. Petrov \cite{petrov:1995}, Theorem
2.3 on p. 51)
yields
\beam\label{eq:ll}
K_{2}
&\le& \P \Big(
\sum_{i=1}^{n\la}  S_{A,i}>x\,(1-\xi)\Big)+
\P\Big(\max_{|m-n \la |\le \delta n}\Big|\sum_{i=m}^{n\la}
S_{A,i}\Big|>x\xi \Big)\nonumber\\
&\le &\P \Big(
\sum_{i=1}^{n\la}  S_{A,i}>x\,(1-\xi)\Big)+
c\,\P \Big(\Big|\sum_{i=1}^{\delta n}
S_{A,i}\Big|>0.5 \xi x \Big)\,.
\eeam
Similarly, using the independence of the  \rv s $(S_{A,i})$ and a
maximal inequality,
\beam\label{eq:ggd}
K_{2}
&\ge &\P \Big(\sum_{i=1}^{\la n}
S_{A,i}>x(1+\xi) \Big)-
c\,\P \Big(\Big|\sum_{i=1}^{\delta n}
S_{A,i}\Big|>0.5 \xi x \Big)\,,
\eeam
where $\delta,\xi$ can be made arbitrarily small provided $n$ is
sufficiently large. Next we give bounds for the \pro ies in
\eqref{eq:ll} and \eqref{eq:ggd}.
We have for any real $s>0$ and $y>0$,
\beao
\P(\sum_{i=1}^{sn }S_{A,i}>y)
&\le &\sum_{i=1}^{sn }
\P\Big(\sum_{k=1}^{s n }S_{A,k}>y,S_{A,i}>y,S_{A,j}\le
  y,j\ne i\Big)\\
&&+
\P(\cup_{k=1,j\ne k}^{s n }\{S_{A,k}>y,S_{A,j}>y\})\\
&\le &s n \P(\sum_{k=1}^{s n }S_{A,k}>y,S_{A,1}>y,S_{A,j}\le
  y,j\ne 1)+
[n
\P_A(S_{A}>y)]^2\\
&\le &s\,n \P_A(S_{A}>y)+
[s\,n
(\P_A(S_{A}>y)]^2\,.
\eeao
Hence, because of the \regvar\ of $X$, uniformly for $x\in\tilde\Lambda_n$,
\beao
\dfrac{\P \Big( \sum_{i=1}^{N_A(n)-1} S_{A,i}>x\Big)}{n\P(|X|>x)}&\le& 
\dfrac{\la \P_A(S_A>x(1-\xi)) +c\delta \P_A (|S_A| >0.5 \xi x) }{\P(|X|>x)}+o\Big(\dfrac{\P_A( S_{A}>x )}{\P(|X|>x)}\Big)\\
&\le &
\dfrac{\lambda\P_A(S_A>x(1-\xi))}{\P(|X|>x)}+o\Big(\dfrac{\P_A( S_{A}>x )}{\P(|X|>x)}\Big)\,.
\eeao
We obtain the last inequality, taking into account that the argument above can be applied to the left tail of $S_A$ as well. This proves the upper bound \eqref{eq:ub}. 

On the other hand, for $s>0$, sufficiently large $n$, small $\vep>0$ 
and $y\in \tilde\Lambda_n$,
\beao
\P(\sum_{i=1}^{sn }S_{A,i}>y)&\ge &
\P\Big(\cup_{i=1}^{sn }\Big\{\sum_{k\neq i}S_{A,k}\le \vep y,S_{A,i}>y(1+\vep),S_{A,j}\le
  y(1+\vep),j\ne i\Big\}\Big)\\
&\ge  & sn \P\Big(\sum_{k=2}^{s n }S_{A,k}\le \vep y,S_{A,1}>y(1+\vep),S_{A,j}\le
  y(1+\vep),j\ne 1\Big)\\
&\ge & (1-\vep) s n  \P_A(S_{A}>y(1+\vep))\,.
\eeao
We conclude from \eqref{eq:ggd} that, uniformly for $x\in\tilde\Lambda_n$,
\beao
\dfrac{\P \Big( \sum_{i=1}^{N_A(n)-1} S_{A,i}>x\Big)}{n\P(|X|>x)}&\ge& 
(1-\vep)\dfrac{\la \P  (
S_{A}>x(1+\xi)(1+\vep) )-c\delta \P_A (|S_A| >0.5 \xi x(1+\vep))}{\P(|X|>x)}\,.\nonumber
\eeao
Now, the lower bound \eqref{eq:lb} is proved in a similar fashion as above.
\end{proof}
In view of Lemmas \ref{lem:1} and \ref{lem:3}, letting first $x\to\infty$ and then $\vep\to0$ and $\xi\to0$ and using \regvar\ of $X$ we obtain 
$$
\dfrac{b_+}{\lambda} = \lim_{x\to\infty }\dfrac{\P_A(S_A>x)}{\P(|X|>x)}\quad\mbox{uniformly on }\tilde\Lambda_n.
$$
In particular this relation holds along the sequences $x_n=cb_n\in\tilde\Lambda_n$ satisfying $x_{n+1}/x_n\to 1$. A sequential version of \regvar\ then implies that $\P_A(S_A>x)$ is \regvary ;
see Bingham et al. \cite{bingham:goldie:teugels:1987}, Theorem 1.9.2.
An application of Theorem \ref{thm:nagaev} and Theorem \ref{thm:mc} finishes the proof of the theorem.
\end{proof}
\bre
Regular variation of $\P_A(S_A>x)$ also implies the following:
\beao
\sup_{x>b_n}\Big|\dfrac{\P\Big(\sum_{i=1}^{N_A(n)-1} 
S_{A,i}>x\Big)}{n\,\P(|X|>x)}- b_+ \Big|\to 0\,.
\eeao 
For the region $x\in\tilde\Lambda_n$ this fact was proved above. Now assume that 
$x\ge \ex^{s_n}$. We have by Theorem~\ref{thm:nagaev} for $\alpha>1$, since $x\ge k$ for $k\le n$, uniformly for $x\ge \ex^{s_n}$,
\beao
\nonumber\P\Big(\sum_{i=1}^{N_A(n)-1} 
S_{A,i}>x\Big)&\sim & \sum_{k=2}^n\P(N_A(n)=k)\,k\,\P_A(S_A>x)\\
\nonumber &\sim &  \P_A(S_A>x)\,\E N_A(n)\sim n\,(\E \tau_A)^{-1}\P_A(S_A>x)\,.
\eeao
An inspection of the proof of Theorem~\ref{thm:int}
now shows why the precise \ld\ principle for $(X_n)$ might in general not hold in the region $(c_n,\infty)$: the first and the last blocks in $S_n$ are always negligible if $\tau_A\le n$. Thus for any $x\ge b_n$ one has
\beam\label{eq:opp}
\dfrac{\P(S_n>x)}{n\P(|X|>x)}\sim b_+ +\dfrac{\P(S_n>x,\tau_A>n)}{n\P(|X|>x)}=b_++r(x).
\eeam
In the region $\Lambda_n$, $r(x)$ is uniformly negligible because it is smaller than $\P(\tau_A>n)/(n\P(|X|>x))$. Therefore the precise large deviation result of Theorem \ref{thm:mc} holds. However, $r(x)$ cannot be neglected in general. It may  influence the very large deviations for $x>c_n$ in a  complicated way:  the Nummelin regeneration scheme cannot be used on $\{\tau_A>n\}$. Below  two special examples of functions of Markov chains are given, where the specific dynamics of the models give some clue on the behavior of the second term.
\ere
\bexam{\rm
Consider the autoregressive process of order 1,
$X_t=\varphi X_{t-1}+ B_t$ for some constant $\varphi\in (-1,1)$ and
an iid \seq\ $(B_t)$
\st\ $B$ is \regvary\ with index $\alpha$ and $\E B=0$ if $\E
|B|<\infty$. It is known from Mikosch and Samorodnitsky
\cite{mikosch:samorodnitsky:2000} that one can choose
$\Lambda_n=(b_n,\infty)$ with $(b_n)$ from Theorem~\ref{thm:nagaev}
and
\beao
b_+= (1-|\varphi|^\alpha)
\Big( \dfrac{p}{(1-\varphi)_+^\alpha} + \dfrac{q}{(1-\varphi)_-^\alpha}\Big)\,,
\eeao
where $p=1-q=\lim_{\xto} \P(B>x)/\P(|B|>x)$. This result was derived
without any further conditions on $B$. The same result follows from
Theorem ~\ref{thm:mc} under more restrictive conditions, e.g.
if $B$ has a non-singular \ds\ with respect to Lebesgue \ms\ (see
Alsmeyer  \cite{alsmeyer:2003}). Thus the remainder term $r(x)$ in \eqref{eq:opp} is uniformly negligible over $(b_n,\infty)$.}
\eexam

\subsection{Solution to stochastic recurrence equations}\label{subsec:sre}
In this section, we consider a special class of stationary
\MC s $(X_t)$ for
which we can apply Theorem~\ref{thm:mc} by considering it as a function of its enlargement $(\Phi_t)$ possessing an atom. Let $((A_t,B_t))_{t\in\bbz}$
be an iid \seq\ \st\ for a generic element $(A,B)$ the following
set of conditions ${\bf SRE }_\alpha$ holds:
\begin{itemize}
\item
$A\ge 0$, $A\ne 0$ a.s., $B\ne 0$ a.s.,
and  the  \ds\ of $(A,B)$ is non-singular with respect to the Lebesgue measure on $\R^2$.
\item The \MC\ $X_t=\Psi_t(X_{t-1})$ is the unique solution to a \sre\  with
iid iterated functions $\Psi_t$ satisfying the
following additional conditions:
\begin{itemize}
\item The Lipschitz coefficients
$L_t$ of the mapping $\Psi_t$ satisfy $\E \log^+ L_t<\infty$.
\item The top Lyapunov exponent of $(\Psi_t)$ is strictly negative.
\item For any $t$,
\beam\label{eq:sre}
A_t\,X_{t-1}- |B_t|\le  X_t\le  A_t\,X_{t-1}+|B_t|\,.
\eeam
\end{itemize}
\item
There exists an $\alpha>0$ \st\ $\E A^\alpha=1$,
  $\E A^{\alpha+\delta}<\infty$ and $\E|B|^{\alpha+\delta}<\infty$ for
  some $\delta>0$.
\item The conditional law of $\log A$, given $A\ne 0$, is
non-arithmetic.
\item
The \ds\ of $X$ is \regvary\ with index $\alpha>0$ in the
following sense:
There exist constants
$c_\infty^+,c_\infty^-\ge 0$ \st\ $c_\infty^++c_\infty^->0$ and
\beam\label{eq:power}
\P (X>x)\sim c_\infty^+\,x^{-\alpha}\,,
\quad\mbox{and}\quad \P(X\le -x)\sim
c_\infty^-\,x^{-\alpha}\,\quad \mbox{as}\; \xto \,.
\eeam
\end{itemize}
These conditions are motivated by the well studied affine case:
\beam\label{eq:sr}
X_t=A_t\,X_{t-1}+B_t\,,\quad t\in\bbz\,.
\eeam
The \sre\ \eqref{eq:sr} has attracted a lot of attention, starting
with pioneering work of Kesten \cite{kesten:1973} who proved that
\eqref{eq:sr} has a stationary solution $(X_t)$ under mild conditions on the
\ds\ of $(A,B)$. This solution has a \regvary\ marginal \ds\ with
index $\alpha>0$ solving the equation $\E A^\kappa=1$, $\kappa>0$.
Kesten's theory was formulated for multivariate $X_t$'s. In the
one-dimensional case, Goldie \cite{goldie:1991} gave an alternative
proof of the \regvar\ of $X$ and he also determined the constants
$c_\infty^-$ and
$c_\infty^+$. In particular, for $B\ge 0$ a.s. he showed that
\beao
c_\infty^+=\dfrac{\E[(B_1+A_1 X_0)^\alpha-(A_1 X_0)^\alpha]}{\alpha  \E A^\alpha\log A}\,.
\eeao
Buraczewski et al.
\cite{buraczewski:damek:mikosch:zienkiewicz:2011} proved a precise \ld\
principle \eqref{eq:1} in the affine case \eqref{eq:sr} in the region
$\Lambda_n=(b_n,c_n)$, where $(b_n)$ is chosen as in Theorem~\ref{thm:mc}
and $c_n=\ex^{s_n}$ for any \seq\ $(s_n)$ \st\ $s_n\to\infty$ and
$s_n=o(n)$. 
The proof in \cite{buraczewski:damek:mikosch:zienkiewicz:2011} is rather technical and uses some
deep analysis of the structure of the random walk $(S_n)$ determined by
the equation \eqref{eq:sr}. In what
follows, we will show that Theorem~\ref{th:main} can be used to
establish the same results  by using the Markov structure of the \seq\
$(X_t)$. The proofs of this section will need less technical
efforts than in \cite{buraczewski:damek:mikosch:zienkiewicz:2011} and give some insight into
precise \ld\ principles for classes of Markov chains larger than
the affine case \eqref{eq:sr}.
\par
Goldie \cite{goldie:1991} already considered \sre s beyond affine
structures. Some of his examples satisfy inequality \eqref{eq:sre}:
\bexam{\rm
Consider the solution to the \sre
\beam\label{eq:gold1}
X_t= \max(A_tX_{t-1},B_t)\,,\quad t \in\bbz.
\eeam
It exists under the conditions $\E \log A<0$, $\E \log^+
B<\infty$ and satisfies \eqref{eq:sre}.
Moreover, if $\E A^\alpha=1$, $\E A^\alpha\log A<\infty$,
the conditional law of $\log A$, given $A\ne 0$, is
non-arithmetic and  $\E(B^+)^\alpha<\infty$, then the unique solution to
\eqref{eq:gold1} satisfies relation \eqref{eq:power}; see Goldie
\cite{goldie:1991},
Theorem~5.2.
}
\eexam
\bexam{\rm
Consider an iid \seq\ $((A_t,C_t,D_t))_{t\in\bbz}$ with a generic
element $(A,C,D)$  \st\ $A\ge 0$ a.s. and $C,D$ are real-valued.
The solution to the equation
\beao
X_t= A_t\,\max(C_t,X_{t-1})+D_t\,,\quad t\in\bbr\,,
\eeao
was considered by Letac \cite{letac:1986}.
It exists under the conditions $\E \log A<0$, $\E \log^+
C<\infty$, $\E \log^+
D<\infty$ and satisfies \eqref{eq:sre} if $D\ge 0$ a.s.
Indeed, if we write $B_t=  A_t \,C_t^++D_t $ then
\beao
|X_t-A_t\,X_{t-1}|\le A_t\,(C_t-X_{t-1})_++D_t\le A_t \,C_t^++D_t=B_t\,.
\eeao
This example is also known to satisfy \eqref{eq:power} (see Goldie
\cite{goldie:1991},
Theorem~6.2): if $A  \ge 0$, $\E (A C^+)^\alpha<\infty$, $\E
|B|^\alpha<\infty$
and $A$ satisfies all conditions of the previous example then
\eqref{eq:power} holds.
}
\eexam
Goldie \cite{goldie:1991} gave various other examples of \sre s
satisfying \eqref{eq:power}. Recently, Mirek \cite{mirek:2011}
considered multivariate analogs of not necessarily affine \sre s
satisfying a condition of type \eqref{eq:sre} (adjusted to the
multivariate case). He proved the \regvar\ of the marginal \ds\
and also gave examples supplementary to those
in  \cite{goldie:1991}. The use of \eqref{eq:sre} in his paper
was also the motivation for us to include in this paper \sre s
which do not necessarily satisfy \eqref{eq:sr}.
\par
In what follows, it will be convenient to write
\beao
\Pi_0=1\quad\mbox{and}\quad \Pi_j=\prod_{i=1}^j A_1\cdots A_j\,,\quad
j\ge 1\,.
\eeao
\bth\label{thm:sre}
Assume that the  stationary \MC\ $(X_t)$ satisfies the condition ${\bf
  SRE}_\alpha$ for some $\alpha>0$ and  $\E X=0$ if $\E |X|<\infty$.
Then the precise \ld\ principle
\eqref{eq:1} holds with
\beam\label{eq:const}
b_+=\E\Big[\Big(1+\sum_{i=1}^\infty \Pi_i \Big)^\alpha-\Big(\sum_{i=1}^\infty \Pi_i\Big)^\alpha\Big]
\eeam
in the regions $\Lambda_n=(b_n,c_n)$ given by
\begin{itemize}
\item
$0<\alpha< 1$: $\Lambda_n=(b_n,  \infty)$ for any $(b_n)$ satisfying
of $b_n/n^{1/\alpha}\to\infty$.
\item
$1< \alpha$ and $\alpha\neq 2$: $\Lambda=(b_n,c_n)$ for any  \seq\ $(b_n)$
satisfying $b_n/n^{1/\alpha\vee 0.5+\delta}\to\infty$
 for any $\delta>0$, and $c_n=\ex^{\gamma n}$ for sufficiently small $\gamma>0$.
\end{itemize}
\ethe
\begin{proof} The condition ${\bf RV}_\alpha$ follows from \regvar\ of
  the marginals. Indeed, iteration of \eqref{eq:sr} yields for fixed
  $d\ge 1$,
\beao
X_0 \Pi_n + R_{n,1} \le X_t\le  X_0 \Pi_n + R_{n,2}\,,\quad n=1,\ldots,d,
\eeao
where $(R_{n,i})_{n=1,\ldots,d}$, $i=1,2,$ is independent of
  $X_0$. Moreover, by the assumptions on $(A,B)$, $\E
  |R_{n,i}|^{\alpha+\delta}<\infty$. Therefore
\beao
\bfX_d=(X_1,\ldots,X_d)= X_0 (\Pi_1,\ldots,\Pi_d)+ \bfR_d\,.
\eeao
Since $X_0$ is assumed \regvary\ with index $\alpha$
an
application
of a multivariate version of a result of Breiman \cite{breiman:1965}
(see Basrak et. al \cite{basrak:davis:mikosch:2002a}) shows that
$ X_0 (\Pi_1,\ldots,\Pi_d)$ is \regvary , and it follows from
Lemma 3.12 in Jessen and Mikosch \cite{jessen:mikosch:2006} and
from $\E|\bfR_d|^{\alpha+\delta}<\infty$ for some $\delta>0$
that $\bfX_d$ is \regvary\ with index $\alpha$. This also means that
one can use the same calculations for $b_+(d)$ given in Bartkiewicz
et al. \cite{bartkiewicz:jakubowski:mikosch:wintenberger:2011}
 and hence the limit $b_+$ exists and is given by the
expression \eqref{eq:const}. Notice that \cite{bartkiewicz:jakubowski:mikosch:wintenberger:2011}
derive the constant $b_+$ only for
$\alpha\in (0,2)$. However, the proofs in the cases
$\alpha\in (1,2)$ and $\alpha>1$ are identical.
\par
Next we verify condition ${\bf AC}_\alpha$ for the region
$(b_n, \infty)$ for any \seq\ $(b_n)$ satisfying
$b_n/n^{1/\alpha}\to \infty$ or, equivalently, $n\,\P(|X|>b_n)\to 0$.
Write $\Pi_{ij}=A_i\cdots A_j$ for any $i,j\in\bbz$ with the convention
that $\Pi_{ij}=1$ if $j,i$.
Iterating \eqref{eq:sre}, we obtain
\beam\label{eq:1a}
X_j\le  \Pi_j\,X_0+\sum_{i=1}^{j}\Pi_{i+1,j}\,|B_{i}|\,,\quad j\ge 0\,.
\eeam
The second term in the \rhs\ of \eqref{eq:1a} is independent of
$X_0$. Hence for $\delta_k>0$,
\beao
\lefteqn{\P(|X_j|>x\delta_k \mid |X_0|>x\delta_k)}\\&\le&
\P(\Pi_j|X_0|>x\delta_k/2\mid
|X_0|>x\delta_k)+\P\Big(\sum_{i=1}^{j}\Pi_{i+1,j}\,|B_{i}|>x\delta_k/2\Big)\\
&=&I_1(x)+I_2(x)\,.
\eeao
Under condition ${\bf SRE}_\alpha$ it follows from Kesten
\cite{kesten:1973} and Goldie \cite{goldie:1991}
that
\beao
Q_j=\sum_{i=-\infty}^{j}\Pi_{i+1,j}\,|B_{i}|<\infty\,,
\eeao
and $(Q_j)$ is the causal solution to the  \sre\ $Q_j= A_j Q_{j-1}+|B_j|$,
$t\in\bbz$, which according to the Kesten-Goldie theory is \regvary\
with index $\alpha$. Therefore
\beao
\sup_{x\ge b_n} n\,I_2(x)= n\,I_2(b_n)\to 0\,,\quad \nto\,,
\eeao
for every $\delta_k>0$ and any \seq\ $(b_n)$ \st\ $b_n /n^{1/\alpha}\to\infty$.
We also have
$$
 I_1(x)\le\frac{\P(\min(\Pi_j,1)\,|X_0|>x\delta_k/2)}{\P(|X_0|>x\delta_k)}.
$$
In view of \eqref{eq:power} there exists a constant $c>0$
such that $\P(X_0>x)\le c\,x^{-\alpha}$, $x>0$.

Using this inequality conditionally on $(A_i)_{1\le i\le j}$, we obtain
$$
\P(\min(\Pi_j,1)\,|X_0|>x\delta_k/2 \mid  (A_i)_{1\le i\le j})\le
c\,(2\min(\Pi_j,1))^\alpha(x\delta_k)^{-\alpha}\,,
$$
and taking expectations,
$$
I_1(x)\le c\,
\E(\min(\Pi_j,1))^\alpha (x\delta_k)^{-\alpha}\,.
$$
Since $\min(y^\alpha,1)\le y^{\alpha-\epsilon}$ for $y\ge 0$,
$\epsilon\in (0,\alpha)$, fixed $\delta_k>0$, and large $n$,
$$
\sup_{x\ge b_n}\delta_k^{-\alpha}\sum_{j= k}^n\P(|X_j|>
x\delta_k\mid |X_0|> x\delta_k)\le c\,
\delta_k^{-2\alpha}\sum_{j= k}^n (\E A^{\alpha-\epsilon})^j.
$$
Since $\E A^{\alpha-\vep}<1$, the \rhs\ is bounded by
$c (\E A^{\alpha-\vep})^k/\delta_k^{2\alpha}$.
Thus ${\bf AC}_\alpha$ is satisfied for any choice of $(b_n)$
with $b_n/n^{1/\alpha}\to \infty$ and $(\delta_k)$ \st\
$(\E A^{\alpha-\vep})^k=o(\delta_k^{2\alpha})$ as $k\to\infty$. In particular, one can
choose $(\delta_k)$ decaying to zero exponentially fast.
\par
Our next goal is to verify \eqref{eq:main}.\\
{\em The case $0<\alpha<1$.} Condition \eqref{eq:main} is immediate from
Remark~\ref{rem:less1}. We can choose $(\delta_k)$ decaying
exponentially fast, as discussed above, and $\vep_k=k^{-2}$.\\[1mm]
{\em The case $\alpha>1$ and $\alpha\neq 2$.} In this case the verification of
\eqref{eq:main} is much more involved. We will employ
Theorem~\ref{thm:mc}. According to this result, we need to verify that
$(X_t)$ is irreducible strongly aperiodic and that
the \MC\ satisfies ${\bf DC}_{p}$ for $p<\alpha$. However, since
$\E A^\alpha=1$, by convexity of the  \fct\ $f(x)=\E A^x$, $x>0$, we have
$f(p)<1$ as $p<\alpha$. Writing $p=\beta k$  where $0<\beta<1$ and $k$ is an integer then 
\beao
\E (|X_1|^{p}-A^p |x|^{p}\mid X_0=x)&\le&\E((A^\beta|x|^\beta+|B|^\beta)^k-(A^\beta |x|^{\beta})^k\mid X_0=x) \\
&=&
\sum_{j=0}^{k-1} {k \choose j}(|x|^\beta)^j\E[(A^\beta)^j (|B|^{\beta})^{k-j}]\\
&\le & c(1+|x|^{p-\beta})\,. 
\eeao
Hence {\bf DC$_p$} is satisfied for any $p<\alpha$.
\par
An application of a result of Alsmeyer \cite{alsmeyer:2003} yields that
the Markov chain $(X_t)$ is aperiodic and  irreducible.
The aperiodicity and $\P$-irreducibility
follow from Theorem 2.1 and Corollary 2.3 in  \cite{alsmeyer:2003}
if and only if the transition kernel of the \MC\
has a component which is absolutely continuous with respect to
Lebesgue \ms . The latter condition is satisfied
in view of the  non-singularity of the \ds\ of $(A,B)$ assumed in
${\bf SRE}_\alpha$ and
since
\beao
\P_x(X>\vep)\ge \P(A\,x-B>\vep)\quad\mbox{and}
\quad \P_x(X\le -\vep)\ge \P(A\,x+B\le -\vep)\quad
\mbox{ for any $\vep>0$.}
\eeao
Thus all assumptions of Theorem~\ref{thm:mc} are satisfied and therefore
its conclusion applies. 
\end{proof}

\subsection{The \garch\ model}
Consider the model \eqref{eq:sv}
with the specification that $(Z_t)$ is an iid symmetric \seq\ and
\beam\label{eq:vol}
\sigma_t^2=\alpha_0+ \sigma_{t-1}^2\,(\alpha_1 Z_{t-1}^2+\beta_1)
=\alpha_0+ \sigma_{t-1}^2 A_t\,,
\eeam
where $\alpha_0,\alpha_1>0$ and $\beta_1\ge 0$. This \sre\ defines
a \garch\ process.  The \garch\
process has been used most frequently for applications in financial
\tsa ; see Andersen et al. \cite{andersen:davis:kreiss:mikosch:2009}.
The theory of Section~\ref{subsec:sre} can be applied to the
affine \sre\ \eqref{eq:vol}.
There exists a unique stationary solution to \eqref{eq:vol}
under the assumption $\E \log A<0$ and
$\sigma$ is \regvary\ under mild conditions on the \ds\ of $Z$.
We will now show a precise \ld\ principle for the process
$(X_t)$
\bth
Consider a \garch\ process $(X_t)$ given by \eqref{eq:sv}
and \eqref{eq:vol} with $\alpha_0,\alpha_1>0$, $\beta_1\in [0,1)$. We assume that there exists an $\alpha>0$, $\alpha\ne 2$ \st:
\begin{itemize}
\item
$Z$ is symmetric with $\var(Z)=1$, $\E |Z|^{\alpha+\delta}<\infty$ for some
$\delta>0$ and the \ds\ of $Z^2$ is non-singular with respect to Lebesgue measure.
\item
There exists an $\alpha>0$ \st\ $\E A^{\alpha/2}=1$.
\end{itemize}
Then the precise \ld\ result \eqref{eq:1} holds in the region  $\Lambda_n=(n^{1/\alpha+\delta},\infty)$ if $\alpha<2$ and $\Lambda_n=(n^{1/2+\delta},\ex^{\gamma n})$ for sufficiently small $\gamma>0$ if $\alpha>2$ with
\beao
b_+=\dfrac{\E[|Z_0+ A_1^{0.5} T_\infty|^{\alpha}-
|A_1^{0.5}T_\infty|^{\alpha}]}{2\E|Z|^{\alpha}}\,,
\eeao
and
$T_\infty=\sum_{t=1}^\infty Z_{t}\,
\prod_{i=1}^{t-1}A_{i+1}^{0.5}\,.
$
\ethe
\begin{proof}
We verify the conditions of Theorem~\ref{th:main}.
Since $(\sigma_t^2)$ satisfies the affine \sre\
\eqref{eq:vol} the assumptions on the \ds\ of $A$ imply that
the conditions of
Goldie \cite{goldie:1991}, Theorem 5.2, are satisfied and therefore
$\sigma$ satisfies the relation $\P(\sigma>x)\sim c_ \infty
x^{-\alpha}$ for some positive $c_\infty$
as $\xto$. Following the argument of the proof on top of p. 366
in  Bartkiewicz et al.
\cite{bartkiewicz:jakubowski:mikosch:wintenberger:2011}, we can show
that for $d\ge 1$,
\beao
\dfrac{\P\Big(\Big|(X_1,\ldots,X_d)-\sigma_0(Z_1 A_1^{0.5},\ldots,Z_d
  \Pi_d^{0.5})\Big|>x\Big)}{\P(|\sigma|>x)}\to 0\,.
\eeao
Observing that $\E |Z_1 A_1^{0.5} |^{\alpha+\delta}<\infty$,
it follows from Lemma 3.12 in Jessen and Mikosch
\cite{jessen:mikosch:2006} and from a generalization of Breiman's
result in Basrak et al. \cite{basrak:davis:mikosch:2002a}
that ${\bf RV}_\alpha$ holds.
\par
The constant  $b_+$ was derived in
\cite{bartkiewicz:jakubowski:mikosch:wintenberger:2011} for $\alpha\in
(0,2)$ but the proof generalizes to arbitrary $\alpha>0$.
\par
As to ${\bf AC}_\alpha$, it follows by the argument leading to
\eqref{eq:1a} that
\beao
\sigma_j^2=\Pi_j \sigma_0^2 + \alpha_0\,\sum_{i=1}^j \Pi_{i+1,j}\,,\quad
j\ge 0\,.
\eeao
Then
\beao
\P(|X_j|>\delta_kx\mid |X_0|>\delta_k x)\le \P\Big(\Pi_j Z_j^2
\sigma_0^2>(\delta_kx)^2/2\mid |X_0|>\delta_k x\Big)
+  \P\Big( Z_j^2\alpha_0\,\sum_{i=1}^j \Pi_{i+1,j}>(\delta_kx)^2/2\Big)\,,
\eeao
and now one can follow the proof of ${\bf AC}_\alpha$ in
Theorem~\ref{thm:sre}. No conditions on $(\delta_k)$ are required so
far and $(b_n)$ is chosen \st\ $b_n/n^{1/\alpha}\to \infty$.
\par
Next we verify \eqref{eq:main}.\\
{\em The case $0<\alpha<2$.} Here one can use
Remark~\ref{rem:garch}.\\[1mm]
{\em The case $\alpha>2$.}
We apply Theorem
\ref{thm:mc} 
to  $ X_t=h(\Phi_t)$, $t\in\bbz$, where the \MC\ $(\Phi_t)$ is an enlargement of the irreducible \MC\ $(X_t,\sigma_t^2)$ possessing an atom $A$.
\end{proof}

\end{document}